# WEAK AND ALMOST SURE LIMITS FOR THE PARABOLIC ANDERSON MODEL WITH HEAVY TAILED POTENTIALS

By Remco van der Hofstad,[1] Peter Mörters[2]
and Nadia Sidorova[3,4]

*Eindhoven University of Technology, University of Bath and
University College London*

We study the parabolic Anderson problem, that is, the heat equation $\partial_t u = \Delta u + \xi u$ on $(0, \infty) \times \mathbb{Z}^d$ with independent identically distributed random potential $\{\xi(z) : z \in \mathbb{Z}^d\}$ and localized initial condition $u(0, x) = \mathbf{1}_0(x)$. Our interest is in the long-term behavior of the random total mass $U(t) = \sum_z u(t, z)$ of the unique nonnegative solution in the case that the distribution of $\xi(0)$ is heavy tailed. For this, we study two paradigm cases of distributions with infinite moment generating functions: the case of polynomial or Pareto tails, and the case of stretched exponential or Weibull tails. In both cases we find asymptotic expansions for the logarithm of the total mass up to the first random term, which we describe in terms of weak limit theorems. In the case of polynomial tails, already the leading term in the expansion is random. For stretched exponential tails, we observe random fluctuations in the almost sure asymptotics of the second term of the expansion, but in the weak sense the fourth term is the first random term of the expansion. The main tool in our proofs is extreme value theory.

Received June 2006; revised July 2007.
[1]Supported in part by the Netherlands Organization for Scientific Research (NWO).
[2]Supported by an Advanced Research Fellowship and by Grant EP/C500229/1 of the Engineering and Physical Sciences Research Council.
[3]Supported by Grant EP/C500229/1 of the EPSRC.
[4]Supported by the European Science Foundation (ESF).

## 1. Introduction.

1.1. *Motivation and background.* We consider the heat equation with random potential on the integer lattice $\mathbb{Z}^d$ and study the unique nonnegative solution to the Cauchy problem with localized initial data:

$$\begin{aligned}
\partial_t u(t,z) &= \Delta^{\mathrm{d}} u(t,z) + \xi(z) u(t,z), & (t,z) &\in (0,\infty) \times \mathbb{Z}^d, \\
u(0,z) &= \mathbf{1}_0(z), & z &\in \mathbb{Z}^d,
\end{aligned}$$
(1.1)

where $\Delta^{\mathrm{d}}$ denotes the discrete Laplacian,

$$(\Delta^{\mathrm{d}} f)(z) = \sum_{y \sim z} [f(y) - f(z)], \qquad z \in \mathbb{Z}^d, f : \mathbb{Z}^d \to \mathbb{R},$$

and the potential $\{\xi(z) : z \in \mathbb{Z}^d\}$ is a collection of independent, identically distributed random variables. The parabolic problem (1.1) is called the *parabolic Anderson model*. It serves as a natural model for random mass transport in a random medium; see, for example, [4] for physical motivation of this problem, [2, 7] for some interesting recent work, and [6] for a recent survey of the field.

A lot of the mathematical interest in the parabolic Anderson model is due to the fact that it is the prime example of a model exhibiting *intermittency*. This means that, for large times $t$, the overwhelming contribution to the total mass

$$U(t) := \sum_{z \in \mathbb{Z}^d} u(t,z)$$

of the solution is concentrated on a subset of $\mathbb{Z}^d$ consisting of a small number of islands located very far from each other. This behavior becomes manifest in the large-time asymptotic behavior of $U(t)$ and of its moments. For example, it has been proposed in the physics literature [13] (see also [6]) to define the model as intermittent if, in our notation, for $p < q$,

$$\lim_{t \uparrow \infty} \frac{(E[U(t)^p])^{1/p}}{(E[U(t)^q])^{1/q}} = 0.$$

The large-time asymptotic behavior of $U(t)$ has been studied in some detail for potentials with finite exponential moments, that is, if $E[\exp\{h\xi(0)\}] < \infty$, for all $h > 0$. Important examples include [1, 3] for the case of bounded potentials, [8, 9] focusing on the vicinity of double-exponential distributions, and [10], which attempts a classification of the potentials according to the long-term behavior of $U(t)$. Most of the existing results approach the problem via the asymptotics of the *moments* of $U(t)$ and *almost sure* results are derived using Borel–Cantelli type arguments.

If the potential fails to have finite exponential moments, then the random variable $U(t)$ fails to have any moments, and new methods have to be found to study its almost sure behavior. It is believed that for such potentials the bulk of the mass $U(t)$ is concentrated in a small number of "extreme" points of the potential. This suggests an approach using extreme value theory. It is this approach to the long-term behavior of the parabolic Anderson model that we follow in this paper.

In all cases of potentials with finite exponential moments, it turns out that the two leading terms in the asymptotic expansion of $\log U(t)$ are *deterministic*, an effect which we did not expect to hold for potentials with heavier tails. Our investigation is motivated by this conjecture, and therefore, we are particularly interested in finding the first (nondegenerate) random term in the asymptotic expansion of $\log U(t)$. For this purpose, we consider two classes of heavy-tailed potentials:

- Potentials with *stretched exponential tail*, or *Weibull potentials*. The distribution function of $\xi(0)$ is given as $F(x) = 1 - e^{-x^\gamma}$ for some positive $\gamma < 1$. This class represents potentials with an intermediately heavy tail.
- Potentials with *polynomial tail*, or *Pareto potentials*. The distribution function of $\xi(0)$ is given as $F(x) = 1 - x^{-\alpha}$, for some $\alpha > d$. This class represents the most heavy-tailed potentials.

Note that the condition $\gamma < 1$ is necessary to make the potentials heavy-tailed, and recall from [8] that the condition $\alpha > d$ is necessary (and sufficient) for the existence of a unique solution of the parabolic Anderson problem.

A fairly complex picture emerges from the main results of this paper, which are formulated precisely in Section 1.2 below:

- In the case of potentials with *polynomial tails*, already the *leading* order term is nondegenerate random, and we determine its asymptotic distribution, if normalised by $t^{\alpha/(\alpha-d)}(\log t)^{-d/(\alpha-d)}$, which is of extremal Fréchet type with shape parameter $\alpha - d$.
- In the case of *stretched exponential tails*, the *first* term in the expansion, which is of order $t(\log t)^{1/\gamma}$, is deterministic. For the *second* term, which is of order $t(\log t)^{1/\gamma - 1} \log \log t$, the almost sure limsup and liminf differ by a constant factor, and the weak limit agrees with the latter. The *third* term in the weak expansion is still deterministic of order $t(\log t)^{1/\gamma - 1} \log \log \log t$. Only the *fourth* term in the weak expansion, which is of order $t(\log t)^{1/\gamma - 1}$, is nondegenerate and properly renormalized converges to a Gumbel distribution.

These results are in line with the underlying belief that for heavy-tailed potentials the bulk of the mass $U(t)$ is concentrated in a small number of "extreme" points of the potential. However, this is not proved here. Attacking this problem requires a wider range of methods and is the subject of ongoing research.

1.2. *Main results of the paper.* We now give precise statements of our results. As we consider two classes of potentials and study two types of convergence for each class, we formulate four theorems. Recall that $U(t)$ is the total mass of the solution of (1.1) and abbreviate

$$L_t := \frac{1}{t} \log U(t).$$

Throughout this paper we denote by $F(x) = P(\xi(0) \leq x)$ the distribution function of $\xi(0)$ and define

$$M_r := \max_{|z| \leq r} \xi(z),$$

where $|\cdot|$ is the 1-norm on $\mathbb{Z}^d$.

The first two theorems are devoted to potentials with polynomial tails.

THEOREM 1.1 (Almost sure asymptotics for Pareto potentials). *Suppose that the distribution of $\xi(0)$ has a polynomial tail, that is, $F(x) = 1 - x^{-\alpha}$, $x \geq 1$, for some $\alpha > d$. Then, almost surely,*

$$\limsup_{t \to \infty} \frac{\log L_t - d/(\alpha - d) \log t}{\log \log t} = -\frac{d-1}{\alpha - d} \qquad \text{for } d > 1,$$

$$\limsup_{t \to \infty} \frac{\log L_t - d/(\alpha - d) \log t}{\log \log \log t} = \frac{1}{\alpha - d} \qquad \text{for } d = 1$$

*and*

$$\liminf_{t \to \infty} \frac{\log L_t - d/(\alpha - d) \log t}{\log \log t} = -\frac{d}{\alpha - d} \qquad \text{for } d \geq 1.$$

Looking at convergence in law, denoted by $\Rightarrow$, we find that the liminf above becomes a limit

$$\frac{\log L_t - d/(\alpha - d) \log t}{\log \log t} \quad \Rightarrow \quad -\frac{d}{\alpha - d} \qquad \text{as } t \uparrow \infty.$$

This follows from a much more precise result, which identifies the order of magnitude of $L_t$ itself and the limit distribution of the rescaled random variable $L_t$.

THEOREM 1.2 (Weak asymptotics for Pareto potentials). *Suppose that the distribution of $\xi(0)$ has a polynomial tail, that is, $F(x) = 1 - x^{-\alpha}$, $x \geq 1$, for some $\alpha > d$. Then, as $t \uparrow \infty$,*

$$L_t \left( \frac{t}{\log t} \right)^{-d/(\alpha - d)} \quad \Rightarrow \quad Y \qquad \text{where } P(Y \leq y) = \exp\{-\theta y^{d-\alpha}\}$$

*and*

$$\theta := \frac{(\alpha - d)^d 2^d B(\alpha - d, d)}{d^d (d-1)!},$$

*where $B$ denotes the beta function.*

REMARK 1. Recall from classical extreme value theory (see, e.g., [5], Table 3.4.2) that the maximum of $t^d$ independent Pareto distributed random variables with shape parameter $\alpha - d$ has qualitatively the same weak asymptotic behavior as our logarithmic total mass $L_t$. An interpretation of this fact is that in the parabolic Anderson model with polynomial potential, the *random fluctuations* of the potential dominate over the *smoothing* effect of the Laplacian.

The results for potentials with polynomial tails prepare the ground for the discussion of the considerably more demanding case of potentials with stretched exponential tails. The next two theorems are the *main results* of this paper.

THEOREM 1.3 (Almost sure asymptotics for Weibull potentials). *Suppose that $\xi(0)$ has distribution function $F(x) = 1 - e^{-x^\gamma}$, $x \geq 0$, for some positive $\gamma < 1$. Then, almost surely,*

$$\limsup_{t \to \infty} \frac{L_t - (d \log t)^{1/\gamma}}{(d \log t)^{1/\gamma - 1} \log \log t} = d(1/\gamma^2 - 1/\gamma) + 1/\gamma,$$

$$\liminf_{t \to \infty} \frac{L_t - (d \log t)^{1/\gamma}}{(d \log t)^{1/\gamma - 1} \log \log t} = d(1/\gamma^2 - 1/\gamma).$$

The difference between liminf and limsup in Theorem 1.3 is due to fluctuations from the liminf-behavior which occur at very rare times. Indeed, we have, as $t \uparrow \infty$,

$$\frac{L_t - (d \log t)^{1/\gamma}}{(d \log t)^{1/\gamma - 1} \log \log t} \;\Rightarrow\; d(1/\gamma^2 - 1/\gamma).$$

This is a consequence of the next theorem, which also extends the expansion in the weak sense up to the first (nondegenerate) random term.

THEOREM 1.4 (Weak asymptotics for Weibull potentials). *Suppose that $\xi(0)$ has distribution function $F(x) = 1 - e^{-x^\gamma}$, $x \geq 0$, for some positive $\gamma < 1$. Then,*

$$(L_t - (d \log t)^{1/\gamma} - d(1/\gamma^2 - 1/\gamma)(d \log t)^{1/\gamma - 1} \log \log t$$
$$+ (d/\gamma)(d \log t)^{1/\gamma - 1} \log \log \log t)((d \log t)^{1/\gamma - 1})^{-1}$$
$$\Rightarrow \; Y,$$

where $Y$ has a Gumbel distribution

$$P(Y \leq y) = \exp\{-\theta e^{-\gamma y}\}$$

with $\theta := 2^d d^{d(1/\gamma - 1)}$.

REMARK 2. The almost sure results of Theorem 1.3 also hold in the case of (standard) exponentially distributed potentials, that is, when $\gamma = 1$. Extending the methods of this paper, Lacoin [11] has shown that in this case

$$\liminf_{t \to \infty} \frac{L_t - d \log t}{\log \log \log t} = -(d+1) \qquad \text{almost surely}$$

and

$$L_t - d \log t + d \log \log \log t \quad \Rightarrow \quad Y,$$

where $Y$ has a Gumbel distribution $P(Y \leq y) = \exp\{-2^d e^{-y+2d}\}$.

1.3. *Outline of the proofs.* Let $(X_s : s \in [0, \infty))$ be the continuous-time simple random walk on $\mathbb{Z}^d$ with generator $\Delta^{\mathrm{d}}$. By $\mathbb{P}_z$ and $\mathbb{E}_z$, we denote the probability measure and the expectation with respect to the walk starting at $z \in \mathbb{Z}^d$. By the Feynman–Kac formula (see, e.g., [8], Theorem 2.1), the unique solution of (1.1) can be expressed as

$$u(t,z) = \mathbb{E}_0\left[\exp\left\{\int_0^t \xi(X_s)\,ds\right\}\mathbf{1}_z(X_t)\right],$$

and the total mass of the solution is hence given by

(1.2) $$U(t) = \mathbb{E}_0\left[\exp\left\{\int_0^t \xi(X_s)\,ds\right\}\right].$$

In this representation, the main contribution to $U(t)$ comes from trajectories of the random walk which, on the one hand, spend a lot of time at sites where the value of the potential is large but, on the other hand, are not too unlikely under the measure $\mathbb{P}_0$. In particular, the contributing trajectories will not visit sites situated too far from the origin. We introduce two variational problems depending on the potential $\xi$:

(1.3) $\quad N(t) := \max_{r>0}\left[M_r - \frac{r}{t}\log\frac{r}{2det}\right] \quad \text{and} \quad \underline{N}(t) := \max_{r>0}\left[M_r - \frac{r}{t}\log M_r\right],$

which reflect the interaction of these two factors. Indeed, up to an additive error which goes to zero, $N(t) - 2d$ is an upper and $\underline{N}(t) - 2d$ a lower bound for $L_t$. For most of our applications these bounds are sufficiently close to each other. Our proofs are based on first making these approximations precise, and then investigating the asymptotics of the random variational problems by means of extreme value theory.

To see the relation between $U(t)$ and the approximating functions in more detail, note that the probability that a continuous-time random walk visits a point $z \in \mathbb{Z}^d$ with $|z| = r \gg \sqrt{t}$ is roughly

$$\mathbb{P}(X_t = z) \lessapprox e^{-2dt}\frac{(2dt)^r}{r!} \approx \exp\left\{-r\log\frac{r}{2det} - 2dt\right\}.$$

If $|X_s| \leq r$ for all $s \in [0,t]$, then $\int_0^t \xi(X_s)\,ds \leq M_r t$. This gives the upper bound

$$L_t = \frac{1}{t}\log U(t) \lessapprox \max_{r>0}\left[M_r - \frac{r}{t}\log\frac{r}{2det}\right] - 2d = N(t) - 2d.$$

For a lower bound, we fix a site $z \in \mathbb{Z}^d$ and $\rho \in (0,1)$, and consider only trajectories which remain constant equal to $z$ during the entire time interval $t(\rho, 1]$. The probability of this strategy is

$$\mathbb{P}(X_s = z \ \forall s \in [\rho t, t]) \gtrapprox \left(\frac{1}{2d}\right)^{|z|} e^{-2d\rho t}\frac{(2d\rho t)^{|z|}}{|z|!}e^{-2d(1-\rho)t}$$

$$\approx \exp\left\{-|z|\log\frac{|z|}{e\rho t} - 2dt\right\},$$

while the contribution of these trajectories to the exponent is $\int_0^t \xi(X_s)\,ds \geq t(1-\rho)\xi(z)$. Optimizing over $z \in \mathbb{Z}^d$ and $\rho \in (0,1)$, we arrive at a lower bound of the form

$$L_t = \frac{1}{t}\log U(t) \gtrapprox \max_{z \in \mathbb{Z}^d}\max_{0<\rho<1}\left[(1-\rho)\xi(z) - \frac{|z|}{t}\log\frac{|z|}{e\rho t}\right] - 2d.$$

Interchanging the maxima over $\rho$ and $z$, and maximizing over $\rho \in (0,1)$, gives

$$\max_{z \in \mathbb{Z}^d}\max_{0<\rho<1}\left[(1-\rho)\xi(z) - \frac{|z|}{t}\log\frac{|z|}{e\rho t}\right] = \max_{|z|<t\xi(z)}\left[\xi(z) - \frac{|z|}{t}\log\xi(z)\right] \approx \underline{N}(t),$$

where the condition $|z| < t\xi(z)$ arises from $\rho < 1$ and can be dropped when $t$ is sufficiently large.

Supposing for the moment that these approximations are sufficiently accurate, we can use extreme value theory to derive asymptotics for $N(t)$ and $\underline{N}(t)$ which then extend to $L_t$. While the almost sure results follow directly from results on the almost sure behavior of maxima of i.i.d. random variables, the key to the weak limit statements is to write $N(t)$ and $\underline{N}(t)$ as a functional of the point process

$$\sum_{z \in \mathbb{Z}^d} \varepsilon_{(z/A_t, (\xi(z)-B_t)/C_t)}$$

for suitable scaling functions $A_t$, $B_t$, $C_t$ and show convergence of the point processes along with the functionals. The nature of our functionals will require a somewhat nonstandard set-up, but the core of the arguments in this part of the proof is using familiar techniques of extreme value theory.

The feasibility of this strategy of proof depends on the quality of the approximation of $L_t$ by $N(t) - 2d$, respectively, $\underline{N}(t) - 2d$. In the case of potentials with polynomial tails, the arguments sketched above show that $L_t/N(t) \to 1$ almost surely, which suffices to infer both weak and almost sure limits of $L_t$ from those of $N(t)$. These arguments are technically less demanding, which allows us to exhibit the strategy of proof very clearly, while in the harder case of potentials with stretched exponential tails technical difficulties may obscure the view to the underlying basic ideas. In the latter case the bounds $\underline{N}(t)$ and $N(t)$ have the same almost sure behavior up to the second term, but their weak behavior when scaled as in Theorem 1.4 differs in that the limiting laws are Gumbel with different location parameter. A considerably refined calculation allows us to show that in probability $L_t$ can be approximated by $\underline{N}(t)$ up to an additive error of order smaller than $(\log t)^{1/\gamma - 1}$, and hence, the weak limit theorem for $L_t$ coincides with that of $N(t)$. This is the most delicate part of the proof, where we rely on a thorough study of the behavior of the potential along random walk paths.

The paper is organized as follows. In Section 2 we prove preliminary results, which will be relevant for both classes of potentials. We start Section 2.1 with Lemma 2.1, where we show that $N(t)$ and $\underline{N}(t)$ are well defined and can be expressed directly in terms of the potential $\xi$. In Lemmas 2.2 and 2.3 we compute upper and lower bounds for $L_t$ in a form which will be simplified to $N(t)$ and $\underline{N}(t)$ in the course of the proofs. In Section 2.2 we prepare the discussion of the extreme value behavior of the bounds with two general lemmas dealing with point processes derived from i.i.d. random variables.

Section 3 is devoted to potentials with polynomial tails. In Section 3.1 we analyze the bounds computed in Section 2.1 and the asymptotic behavior of the optimal value $r$ in the definition of $N(t)$ in (1.3). We infer from this that $L_t/N(t) \to 1$ (see Proposition 3.2) and therefore already the first term of the asymptotic expansion of $L_t$ is nondeterministic. Since $M_r$ is the main ingredient in the definition of $N(t)$, we need to find sharp bounds for $M_r$, which we do in Section 3.2 using extreme value theory. In Section 3.3 we find the weak asymptotics for $N(t)$, and hence of $L_t$, using the point processes technique developed in Section 2.7. This proves Theorem 1.2. Finally, in Section 3.4, we use the bounds for $M_r$ and the weak convergence of $N(t)$ to find the almost sure asymptotics of $N(t)$, and hence of $L_t$, which is stated in Theorem 1.1.

In Section 4 we discuss stretched exponential potentials. This is considerably harder than the polynomial case, and we have to refine the approximation of $L_t$ in several steps. In Section 4.1 we find almost sure bounds for $M_r$ with a high degree of precision, using extreme value theory. Then we show in Proposition 4.2 that $N(t) - 2d$ and $\underline{N}(t) - 2d$ are indeed upper and lower bounds for $L_t$ up to an additive error converging to zero. In Section 4.2 we find weak asymptotics for $\underline{N}(t)$ and $N(t)$, which turn out to be different in the fourth term and are therefore insufficient to give the weak asymptotic for $L_t$. In Section 4.3 we therefore show that $(L_t - \underline{N}(t))/(\log t)^{1/\gamma - 1} \Rightarrow 0$. This approximation, formulated as Proposition 4.6, implies that the weak asymptotics of $\underline{N}(t)$ apply in the same form to $L_t$, and this completes the proof of Theorem 1.4. Finally, in Section 4.4, we study the almost sure behavior of $\underline{N}(t)$ and of $N(t)$, using our knowledge of the behavior of the maximum $M_r$. It turns out that $N(t)$ and $\underline{N}(t)$ are so close to each other that we can get the almost sure upper and lower asymptotics for $L_t$ as stated in Theorem 1.3, additionally using our knowledge of the weak asymptotics from Theorem 1.4.

**2. Notation and preliminary results.** Denote by $\bar{F}(x) := 1 - F(x)$ the tail of the potential and by $J_t$ the number of jumps of the random walk $(X_t : t \geq 0)$ before time $t$. Denote by $\kappa_d(r) r^d$ the number of points in the $d$-dimensional ball of radius $r$ in $\mathbb{Z}^d$ with respect to the 1-norm and $\kappa_d := \lim_{r \to \infty} \kappa_d(r)$. One can easily check that $\kappa_d = 2^d/d!$, but we only need to know that it is nonzero (which follows from the equivalence of all norms on Euclidean space).

Throughout the paper, we use the notation $o(\cdot)$ and $O(\cdot)$ for deterministic functions of one variable (which we specify if there is a risk of confusion). If those functions are allowed to depend on the potential $\xi$ or another variable, then we indicate this by the lower index, writing, for example, $o_\xi$ and $O_\xi$. We say that a family of events $(E_t : t \geq 0)$ holds eventually for all $t$ $\Leftrightarrow$ there exists $T > 0$ such that $E_t$ holds for all $t > T$;

infinitely often $\Leftrightarrow$ there exists a sequence $t_n \uparrow \infty$ such that $E_t$ holds for all $t = t_n$.

2.1. *Bounds for $L_t$ and their properties.* The random functions $N(t)$ and $\underline{N}(t)$ have been defined in terms of the maximum $M_r$. The next lemma provides expressions directly in terms of the potential $\xi$. This, as well as the bounds for $L_t$, which we compute later on, is proved under a mild condition on the growth of the maximum $M_r$ of the potential $\xi$, as $r$ goes to infinity. Later we shall see that this condition is satisfied both for the stretched exponential potentials and for the potentials with polynomial tails.

LEMMA 2.1 [$\underline{N}(t)$ and $N(t)$ in terms of $\xi$]. *Let $\eta \in (0, 1)$. Assume that the distribution of $\xi(0)$ is unbounded from above and that, almost surely, $M_r \leq r^\eta$ eventually for all $r$. Then, almost surely:*

(a) *the maxima $N(t)$ and $\underline{N}(t)$ in* (1.3) *are well defined and the maximizing radii $r(t)$ and $\underline{r}(t)$ satisfy $r(t) \to \infty$ and $\underline{r}(t) \to \infty$ as $t \to \infty$;*

(b) *if $r(t) > 2dt$, then $N(t) = \max_{z \in \mathbb{Z}^d}[\xi(z) - \frac{|z|}{t} \log \frac{|z|}{2det}]$;*

(c) $\underline{N}(t) = \max_{z \in \mathbb{Z}^d}[\xi(z) - \frac{|z|}{t} \log_+ \xi(z)]$ *eventually for all $t$, where $\log_+(x) := \log(x \vee 1)$.*

PROOF. (a) The maxima in $N(t)$ and $\underline{N}(t)$ are attained because $M_r$ is a right continuous step function which grows slower than $\frac{r}{t} \log \frac{r}{2det}$ and $\frac{r}{t} \log M_r$ as $r \to \infty$, for each fixed $t$. Moreover, as the potential distribution is unbounded from above, we have $M_r \to \infty$ as $r \to \infty$. Since $\frac{r}{t} \log \frac{r}{2det} \to 0$ and $\frac{r}{t} \log M_r \to 0$ as $t \to \infty$ for any fixed $r$, we obtain $N(t) \to \infty$ and $\underline{N}(t) \to \infty$. On the other hand, for any $R > 0$ and $t$ large enough, we have

$$\max_{r \leq R}\left[M_r - \frac{r}{t} \log \frac{r}{2det}\right] \leq M_R + \frac{R}{t}\left|\log \frac{R}{2det}\right| \overset{t \to \infty}{\longrightarrow} M_R < \infty$$

and

$$\max_{r \leq R}\left[M_r - \frac{r}{t} \log M_r\right] \leq M_R + \frac{R}{t}|\log M_R| \overset{t \to \infty}{\longrightarrow} M_R < \infty,$$

which implies that $r(t) > R$ and $\underline{r}(t) > R$ eventually.

(b) Observe that at $r = 2dt$ the function $r \mapsto \frac{r}{t} \log \frac{r}{2det}$ takes its minimum, and that it is decreasing on $(0, 2dt)$ and increasing on $(2dt, \infty)$. Denote by $z_t$ a point such that $M_{r(t)} = \xi(z_t)$ and $|z_t| \leq r(t)$. If $|z_t| \leq 2dt$, then $M_{r(t)} = \xi(z_t) = M_{|z_t|} \leq M_{2dt}$ and, hence, by monotonicity of $M_r$, we have $M_{2dt} = M_{r(t)}$. Since, by monotonicity, the value of $\frac{r}{t} \log \frac{r}{2det}$ at $r = r(t) > 2dt$ is strictly greater than its value $-2d$ at $r = 2dt$, we obtain

$$M_{r(t)} - \frac{r(t)}{t} \log \frac{r(t)}{2det} < M_{2dt} + 2d,$$

which is a contradiction to $r(t)$ maximizing $N(t)$ in (1.3). Hence, $2dt < |z_t| \leq r(t)$. Since $M_{r(t)} = \xi(z_t)$, we obtain, again using monotonicity, that

$$M_{r(t)} - \frac{r(t)}{t} \log \frac{r(t)}{2det} \leq \xi(z_t) - \frac{|z_t|}{t} \log \frac{|z_t|}{2det}.$$

This proves the upper bound for $N(t)$. The lower bound is obvious.

(c) Denote by $z_t$ a point such that $\xi(z_t) = M_{\underline{r}(t)}$ and $|z_t| \leq \underline{r}(t)$. Since $\underline{r}(t) \to \infty$, we obtain $M_{\underline{r}(t)} = \xi(z_t) > 1$ eventually. Then, for large $t$,

$$\xi(z_t) - \frac{|z_t|}{t} \log_+ \xi(z_t) \geq \xi(z_t) - \frac{\underline{r}(t)}{t} \log_+ \xi(z_t) = M_{\underline{r}(t)} - \frac{\underline{r}(t)}{t} \log M_{\underline{r}(t)},$$

which proves the upper bound for $\underline{N}(t)$.

To prove the lower bound, note that since $M_r > 1$ eventually and $\frac{r}{t} \log M_r \to 0$ as $t \to \infty$ for any fixed $r$, we obtain $\underline{N}(t) > 1$ eventually. Let us assume $t$ to be large enough so that this is satisfied. Now denote by $z_t$ a point where the maximum of the expression on the right-hand side is taken. If $\xi(z_t) \leq 1$, then

$$\underline{N}(t) > 1 \geq \xi(z_t) = \xi(z_t) - \frac{|z_t|}{t} \log_+ \xi(z_t),$$

which proves the statement. If $\xi(z_t) > 1$, then assume first that $\xi(z_t) \leq |z_t|/t$. Then

$$\xi(z_t) - \frac{|z_t|}{t} \log_+ \xi(z_t) \leq \xi(z_t) - \xi(z_t) \log \xi(z_t) \leq 1 < \underline{N}(t),$$

as $x \mapsto x - x \log x$ is maximal at $x = 1$. Finally, if $\xi(z_t) > |z_t|/t$ and $\xi(z_t) > 1$, we obtain

$$\underline{N}(t) \geq M_{|z_t|} - \frac{|z_t|}{t} \log M_{|z_t|} \geq \xi(z_t) - \frac{|z_t|}{t} \log \xi(z_t),$$

as $x \mapsto x - a \log x$ is increasing on $[a, \infty)$, for any $a > 0$ and, in particular, for $a = |z_t|/t$. $\square$

In the next two lemmas we find almost sure upper and lower bounds for $L_t$, which hold eventually for all $t$.

LEMMA 2.2 (Upper bound for $L_t$). *Let $\eta \in (0,1)$. Assume that almost surely $M_r \leq r^\eta$ eventually for all $r$. Then, almost surely, eventually for all $t$*

$$L_t \leq \max_{r > 0}\left[M_r - \frac{r}{t} \log \frac{r}{2det}\right] - 2d + o(1) = N(t) - 2d + o(1).$$

PROOF. Let $r_0$ be such that $M_r \leq r^\eta$ for all $r > r_0$. Let us fix some $\theta > 1$ and $\beta = (1-\eta)^{-1}(1+\varepsilon)$, $\varepsilon > 0$. Sorting the trajectories of the random walk $X$ according to the number of jumps made before time $t$ and taking into account the facts that $\xi(X_s) \leq M_{J_t}$, for all $s \leq t$, and that $J_t$ has a Poisson distribution with parameter $2dt$, we obtain

$$\mathbb{E}_0\left[\exp\left\{\int_0^t \xi(X_s)\,ds\right\}\right] = \sum_{n=0}^\infty \mathbb{E}_0\left[\exp\left\{\int_0^t \xi(X_s)\,ds\right\}\mathbf{1}_{\{J_t=n\}}\right]$$

(2.1)

$$\leq \sum_{n=0}^\infty e^{tM_n - 2dt}\frac{(2dt)^n}{n!}.$$

We now give an upper bound for the tail of the series on the right. Using Stirling's formula,

(2.2) $$n! = \sqrt{2\pi n}\left(\frac{n}{e}\right)^n e^{\delta(n)}, \qquad \text{with } \lim_{n \uparrow \infty} \delta(n) = 0$$

and the assumption $M_r \leq r^\eta$, we obtain, for all $n > \max\{r_0, t^\beta\}$, that

$$tM_n - 2dt + n\log(2dt) - \log(n!)$$
$$\leq tn^\eta - n\log\frac{n}{2det} - \delta(n)$$
$$\leq tn^\eta\left(1 - \frac{n^{1-\eta}}{t}\log\frac{n}{2det} - \frac{\delta(n)}{tn^\eta}\right)$$
$$\leq tn^\eta\left(1 - t^\varepsilon \log\frac{t^{\beta-1}}{2de} - \frac{\delta(n)}{tn^\eta}\right) \leq -\theta\log n,$$

eventually for all $t$. If $t$ is large enough, then $t^\beta > r_0$ and the last estimate holds for all $n > t^\beta$. Splitting the sum on the right of (2.1) at $n = \lceil t^\beta \rceil$ and noting that $\sum_{n > \lceil t^\beta \rceil} n^{-\theta} = o(1)$, we obtain

$$U(t) = \mathbb{E}_0\left[\exp\left\{\int_0^t \xi(X_s)\,ds\right\}\right] \leq (t^\beta + 1)\max_{0 \leq n \leq t^\beta}\left[\frac{e^{tM_n - 2dt}(2dt)^n}{n!}\right] + o(1)$$

and hence,

$$L_t \leq \max_{0 \leq n \leq t^\beta}\left[M_n - \frac{n}{t}\log\frac{n}{2det} - \frac{1}{t}\log\sqrt{2\pi n} - \frac{\delta(n)}{t}\right] - 2d + o(1)$$
$$= \max_{1 \leq n \leq t^\beta}\left[M_n - \frac{n}{t}\log\frac{n}{2det}\right] - 2d + o(1)$$
$$\leq \max_{r > 0}\left[M_r - \frac{r}{t}\log\frac{r}{2det}\right] - 2d + o(1),$$

which completes the proof of the lemma. □

LEMMA 2.3 (Lower bound for $L_t$). *Let $\eta \in (0,1)$. Assume that almost surely $\xi(0) \geq 0$ and $M_r \leq r^\eta$ eventually for all $r$. Then, almost surely, eventually for all $t$*

(2.3) $$L_t \geq \max_{0 < \rho < 1}\max_{|z| \geq 1}\left[(1-\rho)\xi(z) - \frac{|z|}{t}\log\frac{|z|}{e\rho t}\right] - 2d + o(1).$$

PROOF. Let $r_0$ be such that $M_r \leq r^\eta$ for all $r > r_0$. Let $\rho \in (0,1)$ and $z \in \mathbb{Z}^d$. Denote by

$$A_t^{z,\rho} := \{J_{\rho t} = |z|, X_s = z \;\forall s \in [\rho t, t]\}$$

the event that the random walk $X$ reaches the point $z$ before time $\rho t$, making the minimal possible number of jumps, and stays at $z$ for the rest of the time. Denote by $P_\lambda(\cdot)$ the Poisson distribution with parameter $\lambda$. Considering a smaller event where $X$ reaches $z$ taking some fixed route, we obtain a lower bound on the probability

$$\mathbb{P}_0(A_t^{z,\rho}) \geq \frac{P_{2d\rho t}(|z|)P_{2d(1-\rho)t}(0)}{(2d)^{|z|}} = \frac{e^{-2dt}(2d\rho t)^{|z|}}{(2d)^{|z|}|z|!} = \frac{e^{-2dt}(\rho t)^{|z|}}{|z|!}$$

$$= \exp\left\{-|z|\log\frac{|z|}{e\rho t} - 2dt - \log\sqrt{2\pi|z|} - \delta(|z|)\right\},$$

where $\delta(|z|)$ is taken from Stirling's formula (2.2). As $\xi(z) \geq 0$ almost surely for all $z$, we obtain by the Feynman–Kac formula (1.2), for all $\rho$ and $z$,

$$U(t) = \mathbb{E}_0\left[\exp\left\{\int_0^t \xi(X_s)\,ds\right\}\right] \geq e^{t(1-\rho)\xi(z)}\mathbb{P}_0(A_t^{z,\rho})$$

$$\geq \exp\left\{t(1-\rho)\xi(z) - |z|\log\frac{|z|}{e\rho t} - 2dt - \log\sqrt{2\pi|z|} - \delta(|z|)\right\}.$$

Since $\delta$ is bounded and $\log\sqrt{2\pi|z|} \leq o(t)$ for $|z| \leq t^\beta$ for any fixed positive $\beta$, this implies

$$(2.4) \qquad L_t \geq \max_{0<\rho<1} \max_{1\leq|z|\leq t^\beta}\left[(1-\rho)\xi(z) - \frac{|z|}{t}\log\frac{|z|}{e\rho t}\right] - 2d + o(1).$$

Let $\beta = (1-\eta)^{-1}(1+\varepsilon)$, $\varepsilon > 0$. As $M_r \leq r^\eta$ for all $r > r_0$, for all $|z| > \max\{r_0, t^\beta\}$ we have

$$\max_{|z|>\max\{r_0,t^\beta\}}\left[(1-\rho)\xi(z) - \frac{|z|}{t}\log\frac{|z|}{e\rho t}\right]$$

$$\leq \max_{|z|>\max\{r_0,t^\beta\}}\left[(1-\rho)M_{|z|} - \frac{|z|}{t}\log\frac{|z|}{e\rho t}\right]$$

$$(2.5) \qquad \leq \max_{|z|>\max\{r_0,t^\beta\}}\left[(1-\rho)|z|^\eta - \frac{|z|}{t}\log\frac{|z|}{e\rho t}\right]$$

$$= \max_{|z|>\max\{r_0,t^\beta\}}\left[|z|^\eta\left(1-\rho - \frac{|z|^{1-\eta}}{t}\log\frac{|z|}{e\rho t}\right)\right]$$

$$\leq \max_{|z|>\max\{r_0,t^\beta\}}\left[|z|^\eta\left(1-\rho - t^\varepsilon\log\frac{t^{\beta-1}}{e\rho}\right)\right] < 0,$$

eventually for all $t$. Recall that $L_t \geq 0$ and take $t$ large enough so that $t^\beta > r_0$. Then (2.5) implies that the maximum in (2.4) can be taken over all $z$ instead of $|z| \leq t^\beta$, which implies the statement of the lemma. □

2.2. *Point processes and their transformations.* We employ a point process approach to extreme value theory, using wherever possible the terminology and framework of [12]. In this section we recall the basic setup from [12],

Chapter 3, and add two slightly nonstandard lemmas, that will provide the technique for the proof of our weak convergence results, Theorems 1.2 and 1.4.

We begin with some measure-theoretic notation: A Borel measure on a locally compact space $E$ with countable basis is called *Radon* if it is *locally finite*, that is, all compact sets have finite measure. A Radon measure $\mu$ is called a *point measure* if there exists a finite or countably infinite collection of points $x_1, x_2, \ldots \in E$ such that $\mu = \sum_i \varepsilon_{x_i}$, where $\varepsilon_x$ denotes the Dirac measure at $x$.

A sequence $(\mu_n)$ of Radon measures *converges vaguely* to a Radon measure $\mu$ if $\int f \, d\mu_n \to \int f \, d\mu$ for any continuous function $f \colon E \to [0, \infty)$ with compact support. We denote by $M_p(E)$ the set of point measures on $E$ endowed with the topology of vague convergence. A *point process* is a random element of $M_p(E)$. For any Radon measure $\mu$ there exists a (unique) point process $N$ called the *Poisson process with intensity measure* $\mu$ characterized by the following two properties:

- for any Borel set $A \subset E$, the random variable $N(A)$ is Poisson distributed with parameter $\mu(A)$,
- for any pairwise disjoint Borel sets $A_1, \ldots, A_n$, the random variables $N(A_1), \ldots, N(A_n)$ are independent.

We now suppose that, for any $r > 0$, the random variables $\{X_{r,z}, z \in \mathbb{Z}^d\}$ are independent identically distributed with values in a state space $G$ (again, locally compact with countable basis), and we denote the corresponding probability and expectation by $\mathsf{P}$ and $\mathsf{E}$, respectively. We suppose further that $\mu$ is a Radon measure on $G$ such that

$$r^d \mathsf{P}(X_{r,0} \in \cdot) \xrightarrow{\mathrm{v}} \mu, \tag{2.6}$$

where $\xrightarrow{\mathrm{v}}$ denotes vague convergence. Then we define, for any $r > 0$, a point process $\zeta_r$ on $\mathbb{R}^d \times G$ by

$$\zeta_r = \sum_{z \in \mathbb{Z}^d} \varepsilon_{(z/r, X_{r,z})}. \tag{2.7}$$

Let $\zeta$ be a Poisson process on $\mathbb{R}^d \times G$ with intensity measure $\mathrm{Leb}_d \otimes \mu$, where $\mathrm{Leb}_d$ denotes the Lebesgue measure on $\mathbb{R}^d$. A trivial generalization of the proof of [12], Proposition 3.21 (with $\mathbb{N}$ replaced by $\mathbb{Z}^d$, $[0, \infty)$ by $\mathbb{R}^d$, $\mathrm{Leb}_1$ by $\mathrm{Leb}_d$ and $r \in \mathbb{N}$ by $r > 0$) implies that $\zeta_r$ converges in law to $\zeta$.

Observe that this implies that $\int f \, d\zeta_r$ converges in law to $\int f \, d\zeta$ whenever $f \colon \mathbb{R}^d \times G \to [0, \infty)$ is continuous with compact support. Unfortunately, this is not strong enough for our applications, as we need to consider a class of functions with noncompact support in $\mathbb{R}^d \times G$. To overcome this problem within this framework we look at a *compactification* of the state space.

We let $\dot{\mathbb{R}}^d$ be the one-point compactification of $\mathbb{R}^d$ and work on the space $\dot{\mathbb{R}}^d \times G$. On this space, $\text{Leb}_d \otimes \mu$ is no longer a Radon measure and $\zeta_r$ is no longer a point process, as there are compact sets of infinite measure. However, it turns out that we can define subspaces $H \subset \dot{\mathbb{R}}^d \times G$, for which the convergence result remains true, while the class of compactly supported integrands is sufficiently rich for our applications.

LEMMA 2.4 (Point processes and i.i.d. sequences). *Let $H \subset \dot{\mathbb{R}}^d \times G$ be a locally compact Borel set such that:*

(i) $\text{Leb}_d \otimes \mu|_H$ *is a Radon measure on $H$,*
(ii) *each $\zeta_r|_H$ is a point process in $H$,*
(iii) *the projection of each compact set in $H$ to the second coordinate is compact in $G$,*
(iv) $\lambda_r(dx, dy) := \sum_{z \in \mathbb{Z}^d} \varepsilon_{zr^{-1}}(dx) \mathsf{P}(X_{r,0} \in dy)|_H \xrightarrow{v} \text{Leb}_d(dx) \otimes \mu(dy)|_H.$

*Then $\zeta_r|_H$ converges in law to the Poisson process with intensity measure $\text{Leb}_d \otimes \mu|_H$.*

PROOF. This follows as in the proof of [12], Proposition 3.21, by looking at the Laplace functionals for nonnegative test functions $f$ with compact support in $H$ instead of $\mathbb{R}^d \times G$, using our assumption (iv) in [12], (3.20), and (iii) in [12], (3.21), to obtain the result. □

In our applications we deal with families of transformations of point processes. The next lemma describes the convergence of point processes under families of transformations. Recall that, thanks to our compactification, compact sets may well contain points with infinite components.

LEMMA 2.5 (Transformed point processes). *Let $H \subset \dot{\mathbb{R}}^d \times G$, $H' \subset \dot{\mathbb{R}}^{d+1}$ be locally compact Borel sets and $\eta$ a Radon measure on $H$. Let $\Sigma_t$, for $t > 0$, be a family of point processes in $H$ converging in law to a Poisson process $\Sigma$ on $H$ with intensity measure $\eta$. Finally, let $T$ and $T_t$, for $t > 0$, be measurable mappings from $H$ to $H'$ satisfying the following conditions:*

(i) $T$ *is continuous;*
(ii) *for each compact $K' \subset H'$, there is a compact $K \subset H$ containing $T^{-1}(K')$ and all $T_t^{-1}(K')$;*
(iii) *there exist compact sets $K_n \subset H$ such that $\eta(K_n) \to 0$, and $T_t \to T$ uniformly on each $K_n^c$.*

*Then $\Sigma_t \circ T_t^{-1}$ are point processes in $H'$ converging in law to the Poisson process $\Sigma \circ T^{-1}$ with intensity measure $\eta \circ T^{-1}$.*

PROOF. It follows from (ii) that the preimages of compact sets in $H'$ under $T$ and all $T_t$ are relatively compact in $H$, which implies that $\Sigma_t \circ T_t^{-1}$ are point processes and $\eta \circ T^{-1}$ is a Radon measure. By [12], Proposition 3.7, $\Sigma \circ T^{-1}$ is a Poisson process in $H'$ with intensity measure $\eta \circ T^{-1}$ and so it suffices to prove that $\Sigma_t \circ T_t^{-1} \to \Sigma \circ T^{-1}$ in law, which is equivalent to showing that the corresponding Laplace functionals

$$\Phi_{\Sigma_t \circ T_t^{-1}}(g) := E \exp\left\{-\int g \circ T_t \, d\Sigma_t\right\}$$

converge to $\Phi_{\Sigma \circ T^{-1}}(g)$ for all continuous $g: H' \to [0, \infty)$ with compact support. Note that

(2.8)
$$|\Phi_{\Sigma_t \circ T_t^{-1}}(g) - \Phi_{\Sigma \circ T^{-1}}(g)|$$
$$\leq |\Phi_{\Sigma_t \circ T_t^{-1}}(g) - \Phi_{\Sigma_t \circ T^{-1}}(g)| + |\Phi_{\Sigma_t \circ T^{-1}}(g) - \Phi_{\Sigma \circ T^{-1}}(g)|.$$

Since $T$ satisfies (i) and (ii), it induces a continuous mapping $m \mapsto m \circ T^{-1}$ between the spaces of point measures $M_p(H)$ and $M_p(H')$ (see [12], Proposition 3.18, where the condition $T^{-1}(K')$ being compact for any compact $K'$ can be weakened to relative compactness). Now the weak convergence $\Sigma_t \to \Sigma$ implies the weak convergence $\Sigma_t \circ T^{-1} \to \Sigma \circ T^{-1}$, which is equivalent to the convergence of the Laplace functionals, and so the second term on the right-hand side of (2.8) converges to zero.

Let us prove that the first term also does. According to (ii), denote by $K \subset H$ a compact set containing $T^{-1}(\operatorname{supp} g)$ and all $T_t^{-1}(\operatorname{supp} g)$. Further, denote $A_t(a, n) = \{\Sigma_t(K_n \cap K) = 0, \Sigma_t(K_n^c \cap K) \leq a\}$. For any fixed $\varepsilon > 0$, one can fix $n$ and $a$ large enough so that

$$\mathsf{P}(A_t(a, n)) \geq \mathsf{P}(\Sigma_t(K_n) = 0, \Sigma_t(K) \leq a)$$
$$\stackrel{t \to \infty}{\longrightarrow} \mathsf{P}(\Sigma(K_n) = 0, \Sigma(K) \leq a) > 1 - \varepsilon,$$

since $\eta(K_n) \to 0$. This implies $\mathsf{P}(A_t(a, n)) > 1 - \varepsilon$, eventually for all $t$. Further, since $T_t \to T$ uniformly on $K_n^c$ and $g$ is compactly supported and so uniformly continuous, we have

$$\sup_{x \in K_n^c} |(g \circ T_t)(x) - (g \circ T)(x)| \leq \varepsilon/a,$$

eventually for all $t$. Hence, we obtain, also using that $g \geq 0$,

$$|\Phi_{\Sigma_t \circ T_t^{-1}}(g) - \Phi_{\Sigma_t \circ T^{-1}}(g)|$$
$$= |\Phi_{\Sigma_t}(g \circ T_t) - \Phi_{\Sigma_t}(g \circ T)|$$
$$\leq 2\mathsf{P}(A_t(a, n)^c) + \mathsf{E}[|e^{-\int_H (g \circ T_t) \, d\Sigma_t} - e^{-\int_H (g \circ T) \, d\Sigma_t}| \mathbf{1}_{A_t(a,n)}]$$
$$\leq 2\varepsilon + \mathsf{E}[(e^{\int_H |(g \circ T) - (g \circ T_t)| \, d\Sigma_t} - 1) e^{-\int_H (g \circ T) \, d\Sigma_t} \mathbf{1}_{A_t(a,n)}]$$
$$\leq 2\varepsilon + (e^\varepsilon - 1) < 4\varepsilon,$$

eventually for all $t$, if $\varepsilon$ is small enough. $\square$

**3. Potentials with polynomial tails.** In this section we consider Pareto potentials, that is, we assume that the distribution function $F$ of $\xi(0)$ is given by $F(x) = 1 - x^{-\alpha}$ for $x \geq 1$, for some $\alpha > d$. In Section 3.1 we show that the upper and lower bounds on $L_t$, which we found in Section 2, are equivalent, so that it suffices to consider the behavior of one of the bounds. Since both bounds are given in terms of the maxima $M_r$, we compute some bounds on $M_r$ in Section 3.2. Then we study the almost sure behavior of the upper bound, $N(t)$, in Section 3.4, and its weak asymptotics in Section 3.3.

3.1. *Asymptotic equivalence of $L_t$ and $N(t)$.* We first investigate the growth of any radius $r(t)$ maximizing the variational problem $N(t)$ in (1.3).

LEMMA 3.1 [Rough asymptotics of $r(t)$ and $M_{r(t)}$]. *Almost surely, as $t \to \infty$:*

(a) $r(t) = t^{\alpha/(\alpha-d)+o_\xi(1)}$,
(b) $M_{r(t)} = (\frac{r(t)}{t})^{1+o_\xi(1)}$.

PROOF. By [8], Lemma 4.2 (with no assumption on the distribution), we have that almost surely

$$-\log[1 - F(M_r)] = d[\log r](1 + o_\xi(1)) \qquad \text{as } r \to \infty,$$

which implies $M_r = r^{d/\alpha + o_\xi(1)}$. Let $r(t)$ be a maximizer of

$$\Phi_t(r) := M_r - \frac{r}{t} \log \frac{r}{2det}.$$

To prove the lower bound on $r(t)$ by contradiction, we assume that $r(t) < t^{(1-\varepsilon)\alpha/(\alpha-d)}$, for fixed $\varepsilon > 0$ and some arbitrarily large $t$. Consider $\bar{r}(t) := t^{(1-\bar{\varepsilon})\alpha/(\alpha-d)}$, for some $\bar{\varepsilon} < \varepsilon \wedge (d/\alpha)$. Using $\frac{r}{t} \log \frac{r}{2det} \geq -2d$, $r(t) \to \infty$ and the asymptotics of $M_r$ for $r \to \infty$,

$$\Phi_t(\bar{r}(t)) - \Phi_t(r(t)) \geq \bar{r}(t)^{d/\alpha-\delta} - r(t)^{d/\alpha+\delta} - \frac{\bar{r}(t)}{t} \log \frac{\bar{r}(t)}{2det} + \frac{r(t)}{t} \log \frac{r(t)}{2det}$$

$$\geq t^{(1-\bar{\varepsilon})\alpha/(\alpha-d)(d/\alpha-\delta)} - t^{(1-\varepsilon)\alpha/(\alpha-d)(d/\alpha+\delta)}$$

$$+ t^{(1-\bar{\varepsilon})\alpha/(\alpha-d)-1} \log[t^{(1-\bar{\varepsilon})\alpha/(\alpha-d)-1}(2de)^{-1}] - 2d > 0,$$

for sufficiently large $t$, if we choose $\delta$ in such a way that the first term dominates. Hence, under our assumption $r(t)$ is not maximizing, which is a contradiction.

To obtain the upper bound on $r(t)$, we assume that $r(t) > t^{(1+\varepsilon)\alpha/(\alpha-d)}$ for large $t$ and again derive a contradiction. Take $0 < \delta < (\varepsilon(\alpha-d))/(\alpha(1+\varepsilon))$.

By the asymptotics of $M_r$ and the fact that $r \mapsto (r^{1-d/\alpha-\delta}/t) \log \frac{r}{2det}$ is decreasing in $r$ for $r \geq t^{(1+\varepsilon)\alpha/(\alpha-d)}$, we obtain

$$\Phi_t(r(t)) \leq r(t)^{d/\alpha+\delta} - \frac{r(t)}{t} \log \frac{r(t)}{2det} = r(t)^{d/\alpha+\delta} \left(1 - \frac{r(t)^{1-d/\alpha-\delta}}{t} \log \frac{r(t)}{2det}\right)$$

$$< r(t)^{d/\alpha+\delta} (1 - t^{\varepsilon-\delta\alpha/(\alpha-d)(1+\varepsilon)} \log[t^{\varepsilon-\delta\alpha/(\alpha-d)(1+\varepsilon)}(2de)^{-1}]) < 0,$$

for sufficiently large $t$ as the power of $t$ is positive by the choice of $\delta$. Hence, $r(t)$ cannot maximize $N(t)$, so that we again obtain a contradiction, thus proving (a).

To prove the asymptotics for $M_{r(t)}$, recall from the beginning of the proof that $M_r = r^{d/\alpha + o_\xi(1)}$. Using also (a), we get, as $t \to \infty$,

$$M_{r(t)} = M_{t^{\alpha/(\alpha-d)+o_\xi(1)}} = t^{d/(\alpha-d)+o_\xi(1)} = \left(\frac{r(t)}{t}\right)^{1+o_\xi(1)},$$

which completes the proof. □

PROPOSITION 3.2 [Almost sure equivalence of $L_t$ and $N(t)$]. *Almost surely,*

$$L_t = (1 + o_\xi(1))N(t).$$

PROOF. By Lemma 3.1 and recalling that $d < \alpha$, we see that the conditions in Lemmas 2.2 and 2.3 are fulfilled. The former lemma gives the upper bound $L_t \leq N(t) + O(1)$ eventually for all $t$. The lower bound can be obtained from the latter lemma in the following way.

By definition of $r(t)$, we have $M_{r(t)} - \frac{r(t)}{t} \log \frac{r(t)}{2det} > 0$, and hence, the second term grows no faster than (but possibly at the same rate as) the first term, that is,

$$(3.1) \qquad 0 \leq \limsup_{t \to \infty} \frac{1}{M_{r(t)}} \frac{r(t)}{t} \log \frac{r(t)}{2det} \leq 1.$$

Let $\varepsilon_t = o_\xi(1)$ be such that $\log \varepsilon_t = o_\xi(\log \frac{r(t)}{t})$, for example, $\varepsilon_t = (\log t)^{-1}$, and define

$$\rho_t = \frac{\varepsilon_t}{M_{r(t)}} \frac{r(t)}{t} \log \frac{r(t)}{2det},$$

which goes to zero according to (3.1). Now, by Lemma 2.3, we obtain

$$L_t \geq (1 - \rho_t)\xi(z_t) - \frac{|z_t|}{t} \log \frac{|z_t|}{e\rho_t t} + O(1)$$

$$= M_{r(t)} - \frac{r(t)}{t} \log \frac{r(t)}{2det} - \rho_t M_{r(t)} + \frac{r(t)}{t} \log \frac{e\rho_t}{2d} + O(1).$$

Since $\rho_t = o_\xi(1)$, the third term is dominated by the first two terms. Also, the fourth one is dominated by the first two terms, since, by Lemma 3.1(b), we have

$$\log \frac{e\rho_t}{2d} = \log \varepsilon_t + (1 + o_\xi(1)) \log \frac{r(t)}{t} - \log M_{r(t)} = o_\xi(1) \log \frac{r(t)}{t},$$

since $\log \varepsilon_t = o_\xi(\log \frac{r(t)}{t})$. Hence, we obtain

$$L_t \geq \left[ M_{r(t)} - \frac{r(t)}{t} \log \frac{r(t)}{2det} \right](1 + o_\xi(1)) = N(t)(1 + o_\xi(1))$$

eventually, which completes the proof. □

3.2. *Bounds on the maximum $M_r$.* Since $L_t$ is controlled by $N(t)$, which is defined in terms of the maxima $M_r$, we need quite sharp almost sure upper and lower bounds for $M_r$. These are derived in this section using a standard technique for independent identically distributed sequences.

Let $(X_i : i \in \mathbb{N})$ under $\mathsf{P}$ be a family of independent identically distributed random variables with distribution function $F(x) = 1 - x^{-\alpha}$, $x \geq 1$. Define the maximum process $(\bar{X}_n : n \in \mathbb{N})$ by $\bar{X}_n = \max_{i \leq n} X_i$. For any $\rho \in \mathbb{R}$ and $c > 0$, let

$$u_\rho(n) = n^{1/\alpha}(\log n)^{1/\alpha}(\log \log n)^{1/\alpha}(\log \log \log n)^{1/\alpha + \rho},$$

$$v_c(n) = cn^{1/\alpha}(\log \log n)^{-1/\alpha}.$$

The next lemma gives bounds for the maximum process in terms of the sequences $(u_\rho(n))$ and $(v_c(n))$ with suitable parameters $\rho$ and $c$.

LEMMA 3.3 (Extreme value theory for Pareto variables). *Almost surely, as $n \to \infty$:*

(i) $\bar{X}_n \leq u_\rho(n)$ *eventually, if $\rho > 0$;*
(ii) $\bar{X}_n \geq u_\rho(n)$ *infinitely often, if $\rho \leq 0$;*
(iii) $\bar{X}_n \geq v_c(n)$ *eventually, if $c < 1$;*
(iv) $\bar{X}_n \leq v_c(n)$ *infinitely often, if $c \geq 1$.*

PROOF. By [5], Theorem 3.5.1, we have the following criterion for the sequence $u_\rho(n)$ to be an eventual upper bound for $\bar{X}_n$:

$$\mathsf{P}(\bar{X}_n \leq u_\rho(n) \text{ ev.}) = 1 \iff \sum_{n=1}^\infty \mathsf{P}(X_1 > u_\rho(n)) < \infty,$$

$$\mathsf{P}(\bar{X}_n \geq u_\rho(n) \text{ i.o.}) = 1 \iff \sum_{n=1}^\infty \mathsf{P}(X_1 > u_\rho(n)) = \infty,$$

where, for a sequence of events $(E_n : n \geq 0)$, we abbreviate $\{E_n \text{ ev.}\}$ for the event that the events of the sequence occur eventually for all $n$, and $\{E_n \text{ i.o.}\}$ if they occur infinitely often. Compute

$$\sum_{n=1}^{\infty} \mathsf{P}(X_1 > u_\rho(n)) = \sum_{n=1}^{\infty} u_\rho(n)^{-\alpha} = \sum_{n=1}^{\infty} \frac{1}{n(\log n)(\log \log n)(\log \log \log n)^{1+\alpha\rho}}.$$

If $\rho > 0$, then the series converges and, hence, $\mathsf{P}(\bar{X}_n \leq u_\rho(n) \text{ ev.}) = 1$. For $\rho \leq 0$, the series diverges and so $\mathsf{P}(\bar{X}_n \geq u_\rho(n) \text{ i.o.}) = 1$.

To prove the last two bounds, recall that $\bar{F} = 1 - F$ denotes the tail of the distribution and note that $\bar{F}(v_c(n)) \to 0$ as $v_c(n) \to \infty$, and

$$\bar{F}(v_c(n)) = v_c(n)^{-\alpha} = c^{-\alpha} n^{-1} \log \log n.$$

This implies $n\bar{F}(v_c(n)) = c^{-\alpha} \log \log n \to \infty$, and hence, we can apply the criterion [5], Theorem 3.5.2, for the sequence $v_c(n)$ to be an eventual lower bound, which says that

$$\mathsf{P}(\bar{X}_n \leq v_c(n) \text{ i.o.}) = 1 \iff \sum_{n=1}^{\infty} \bar{F}(v_c(n)) \exp\{-n\bar{F}(v_c(n))\} = \infty,$$

$$\mathsf{P}(\bar{X}_n \geq v_c(n) \text{ ev.}) = 1 \iff \sum_{n=1}^{\infty} \bar{F}(v_c(n)) \exp\{-n\bar{F}(v_c(n))\} < \infty.$$

We find that

$$\sum_{n=1}^{\infty} \bar{F}(v_c(n)) \exp\{-n\bar{F}(v_c(n))\} = c^{-\alpha} \sum_{n=1}^{\infty} \frac{\log \log n}{n \exp\{c^{-\alpha} \log \log n\}}$$

$$= c^{-\alpha} \sum_{n=1}^{\infty} \frac{\log \log n}{n(\log n)^{c^{-\alpha}}}.$$

This series is finite if and only if $c^{-\alpha} > 1$, which is the case if $c < 1$. Then $\mathsf{P}(\bar{X}_n \geq v_c(n) \text{ ev.}) = 1$. If $c \geq 1$, on the other hand, then $c^{-\alpha} \leq 1$ and the series diverges, which gives $\mathsf{P}(\bar{X}_n \leq v_c(n) \text{ i.o.}) = 1$. □

We would like to use Lemma 3.3 to estimate $M_r = \max_{|z| \leq r} \xi(z)$, which requires looking at a certain subsequence of the maximum process. We can do this easily for the eventual estimates, but some more work is needed to show that events occurring infinitely often will occur infinitely often along the subsequence.

LEMMA 3.4 (Extreme values along sparse subsequences). *Let $(n_k)$ be an increasing sequence of natural numbers such that, for some $\eta \in (0,1)$, $n_{k+1} - n_k < n_k^\eta$ for all $k$ large enough. Then almost surely, as $k \to \infty$:*

(i) $\bar{X}_{n_k} \geq u_\rho(n_k)$ *infinitely often if $\rho < 0$;*

(ii) $\bar{X}_{n_k} \leq v_c(n_k)$ *infinitely often if* $c > 1$.

PROOF. Take $\rho' \in (\rho, 0)$. By the second statement of Lemma 3.3, we have
$$u_\rho(n + n^\eta) = u_\rho(n)(1 + o(1)) \leq u_{\rho'}(n) \leq \bar{X}_n,$$
where the first inequality holds eventually, and the second one infinitely often. Hence, there exists an arbitrarily large $n$ such that $\bar{X}_n \geq u_\rho(n + n^\eta)$. Choose $k$ in such a way that $n_{k-1} < n \leq n_k$. If $n$, respectively, $k$, is large enough, then $n_k - n < n_k - n_{k-1} < n_{k-1}^\eta < n^\eta$. By monotonicity, we obtain $\bar{X}_{n_k} \geq \bar{X}_n \geq u_\rho(n + n^\eta) \geq u_\rho(n_k)$, which proves (i). The argument for (ii) is analogous. □

Now we are able to prove the estimates for $M_r$.

LEMMA 3.5 (Bounds on $M_r$). *Let $\delta > 0$. There are $c_1, c_2 > 0$ such that, almost surely, as $r \to \infty$, the following estimates hold:*

$M_r \leq r^{d/\alpha}(\log r)^{1/\alpha}(\log \log r)^{1/\alpha}(\log \log \log r)^{1/\alpha + \delta}$ *eventually for all $r$,*

$M_r \geq r^{d/\alpha}(\log r)^{1/\alpha}(\log \log r)^{1/\alpha}(\log \log \log r)^{1/\alpha - \delta}$ *for infinitely many $r$,*

$M_r \geq c_1 r^{d/\alpha}(\log \log r)^{-1/\alpha}$ *eventually for all $r$,*

$M_r \leq c_2 r^{d/\alpha}(\log \log r)^{-1/\alpha}$ *for infinitely many $r$.*

PROOF. The eventual inequalities follow directly from Lemma 3.3. Indeed, let $\rho \in (0, \delta)$. As the ball of radius $r$ contains $\kappa_d(r)r^d$ points, we obtain
$$M_r \leq u_\rho(\kappa_d(r)r^d) \leq r^{d/\alpha}(\log r)^{1/\alpha}(\log \log r)^{1/\alpha}(\log \log \log r)^{1/\alpha + \delta},$$
eventually for all $r$. Similarly, taking $c < 1$, we obtain, for a suitable constant $c_1 > 0$,
$$M_r \geq v_c(\kappa_d(r)r^d) \geq c_1 r^{d/\alpha}(\log \log r)^{-1/\alpha} \quad \text{eventually for all } r.$$
To prove the second and the fourth inequality, we need to check that the subsequence $n_r = \kappa_d(r)r^d$ satisfies the condition in Lemma 3.4. Indeed, for $\eta \in (\frac{d-1}{d}, 1)$, we have
$$\lim_{r \to \infty} \frac{n_{r+1} - n_r}{n_r^\eta} = \lim_{r \to \infty} \frac{\kappa_d(r+1)(r+1)^d - \kappa_d(r)r^d}{\kappa_d(r)^\eta r^{d\eta}} = \kappa_d^{1-\eta} \lim_{r \to \infty} r^{d-1-d\eta} = 0.$$
Taking $\rho \in (-\delta, 0)$, we obtain, by Lemma 3.4, for infinitely many $r$
$$M_r \geq u_\rho(\kappa_d(r)r^d) \geq r^{d/\alpha}(\log r)^{1/\alpha}(\log \log r)^{1/\alpha}(\log \log \log r)^{1/\alpha - \delta},$$
as the first inequality holds infinitely often and the last one eventually. Further, taking $c \geq 1$ and a suitable constant $c_2$, we get for the fourth estimate
$$M_r \leq v_c(\kappa_d(r)r^d) \leq c_2 r^{d/\alpha}(\log \log r)^{-1/\alpha} \quad \text{for infinitely many } r,$$
as the first inequality holds infinitely often and the last one eventually. □

3.3. *Weak convergence of $N(t)$*. We now study the weak asymptotics for $N(t)$, which by Proposition 3.2 agrees with those of $L_t$. Our main tool is the point processes technique and we use the preliminary lemmas proved in Section 2.2. Denote

$$(3.2) \qquad a_t := \left(\frac{t}{\log t}\right)^{d/(\alpha-d)} \quad \text{and} \quad r_t := \left(\frac{t}{\log t}\right)^{\alpha/(\alpha-d)},$$

which turn out to be the right scaling for $N(t)$ and for the radius $r(t)$, where the maximum $N(t)$ is attained. Let $G := (0, \infty]$,

$$X_{r,z} := \frac{\xi(z)}{r^{d/\alpha}} \quad \text{and} \quad \mu(dy) := \frac{\alpha\, dy}{y^{\alpha+1}},$$

which is a Radon measure on $G$. For $x > 0$, we have

$$(3.3) \qquad r^d P(X_{r,0} \geq x) = r^d P(\xi(0) \geq x r^{d/\alpha}) = x^{-\alpha} = \mu([x, \infty])$$

and so condition (2.6) is satisfied. Define $\zeta_r$ by (2.7), let $q := d/(\alpha - d)$ and define a locally compact Borel set

$$H := \{(x, y) \in \dot{\mathbb{R}}^d \times G : y \geq q|x|/2\} \cup \{(\infty, \infty)\}.$$

Denote by $B_\varepsilon := \{(x, y) \in \mathbb{R}^d \times G : |x| < \varepsilon, |y| < q\varepsilon/2\}$ a collection of neighborhoods of zero and by $K_n := \{(x, y) \in H : |x| \geq n\}$ a collection of compact sets in $H$.

LEMMA 3.6 (Convergence of the point processes). $\zeta_r|_H \Rightarrow \zeta|_H$, *where* $\zeta|_H$ *denotes a Poisson process on $H$ with intensity measure* $\eta := \mathrm{Leb}_d \otimes \mu|_H$.

PROOF. It suffices to show that the conditions (i)–(iv) of Lemma 2.4 are satisfied.

(i) To show that $\eta$ is a Radon measure on $H$, it suffices to see that it is finite on the complements of $B_\varepsilon$. This follows as, using $\alpha > d$,

$$\eta(H \setminus B_\varepsilon) = \int_{\mathbb{R}^d} dx \int_{q(\varepsilon \vee |x|)/2}^\infty \frac{\alpha\, dy}{y^{\alpha+1}} = (2/q)^\alpha \int_{\mathbb{R}^d} \frac{dx}{(\varepsilon \vee |x|)^\alpha} < \infty.$$

(ii) To prove that $\zeta_r|_H$ is a point process in $H$, it suffices to show that $\zeta_r(H) < \infty$ almost surely. This follows from the fact that, by Lemma 3.5, almost surely $\xi(z) \leq M_{|z|} < |z|^{d/\alpha+\delta}$ and, hence, $\xi(z) r^{-d/\alpha} < q|z|/(2r)$ for all $|z|$ large enough, so that $(|z|/r, X_{r,z}) \in H$ for only finitely many $z$.

(iii) Obviously, the projection of a compact set in $H$ to the first coordinate is compact in $G$.

(iv) The convergence $\lambda_r \xrightarrow{v} \eta$ with respect to continuous test functions $f \colon \mathbb{R}^d \times G \to [0, \infty)$ with compact support in $\mathbb{R}^d \times G$ follows from (3.3). To prove the vague convergence on $H$, it suffices to show additionally that for each $\varepsilon > 0$ there is $n$ such that $\lambda_r(K_n) < \varepsilon$ eventually for all $r$. We have

$$\lambda_r(K_n) = \sum_{|z| \geq nr} P\left(X_{r,0} \geq q\frac{|z|}{2r}\right)$$

$$= \left(\frac{2}{q}\right)^\alpha \sum_{|z| \geq nr} \left|\frac{z}{r}\right|^{-\alpha} r^{-d} \stackrel{r \to \infty}{\Longrightarrow} \left(\frac{2}{q}\right)^\alpha \int_{\{|x| \geq n\}} |x|^{-\alpha} dx.$$

As $\alpha > d$, the integral converges and the right-hand side is smaller than $< \varepsilon$ if $n$ is large. $\square$

Recall from (3.2) that $a_t = (r_t)^{d/\alpha}$. Let

$$\Psi_t(z) := \xi(z) - \frac{|z|}{t} \log \frac{|z|}{2det}$$

and define the corresponding point process by

$$\Pi_t := \sum_{\{z \in \mathbb{Z}^d \colon \Psi_t(z) > 0\}} \varepsilon_{(z/r_t, \Psi_t(z)/a_t)}.$$

LEMMA 3.7 (Convergence of the transformed point processes). *For each $t$, $\Pi_t$ is a point process on*

$$\hat{H} := \dot{\mathbb{R}}^{d+1} \setminus ((\mathbb{R}^d \times (-\infty, 0)) \cup \{(0,0)\}).$$

*As $t \to \infty$, $\Pi_t$ converges in law to a Poisson process $\Pi$ on $\hat{H}$ with intensity measure*

$$\nu(dx, dy) = dx \otimes \frac{\alpha}{(y + q|x|)^{\alpha+1}} \mathbf{1}_{\{y > 0\}} dy.$$

PROOF. Observe that

$$\frac{\Psi_t(z)}{a_t} = \frac{\xi(z)}{a_t} - \frac{|z|}{ta_t} \log \frac{|z|}{2det} = \frac{\xi(z)}{a_t} - (q + o(1))\left|\frac{z}{r_t}\right| - \frac{|z/r_t| \log |z/r_t|}{\log t}$$

and hence,

$$\Pi_t = (\zeta_{r_t}|_H \circ T_t^{-1})|_{\hat{H}} \qquad \text{eventually for all } t,$$

for a transformation $T_t \colon H \to H' := \dot{\mathbb{R}}^{d+1} \setminus \{0\}$ given by

$$T_t \colon (x, y) \mapsto \begin{cases} (x, y - q|x| - \delta(t, x)), & \text{if } x \neq \infty \text{ and } y \neq \infty, \\ \infty, & \text{otherwise,} \end{cases}$$

where $\delta(t,x) \to 0$ as $t \to \infty$ uniformly on all $K_n^c$. Finally, we define $T: H \to H'$ by $T(x,y) = (x, y - q|x|)$ if $x \neq \infty$ and $y \neq \infty$ and $T(x,y) = \infty$ otherwise. It now suffices to show that

$$\zeta_{r_t}|_H \circ T_t^{-1} \quad \Longrightarrow \quad \zeta|_H \circ T^{-1},$$

as the Poisson process on the right has the required intensity by a straightforward change of coordinates. This convergence holds by Lemmas 3.6 and 2.5, provided the conditions (i)–(iii) of the latter are satisfied, which we now check:

(i) $T$ is obviously continuous.

(ii) For each compact set $K' \subset H'$, there is an open neighborhood $V'$ of zero such that $K' \subset H' \setminus V'$. Since $T_t \to T$ uniformly in $K_n^c$ and since $T(K_n) \cap V' = T_t(K_n) = \varnothing$ for large $n$, there exists an open neighborhood $V \subset H$ of zero such that $T(V) \subset V'$ and $T_t(V) \subset V'$ for all $t$. Hence, for $K := H \setminus V$, we obtain $T^{-1}(K') \subset T^{-1}(H' \setminus V') \subset K$ and, similarly, $T_t^{-1}(K') \subset K$ for all $t$.

(iii) Recall that $\delta(x,t) \to 0$ uniformly on $K_n^c$ and that $\eta(K_n) \to 0$, as $\eta$ is finite away from zero. $\square$

The following proposition, together with Proposition 3.2, completes the proof of Theorem 1.2.

PROPOSITION 3.8 [Weak asymptotics for $N(t)$]. With $a_t$ as in (3.2) and $\theta$ as in Theorem 1.2,

$$\frac{N(t)}{a_t} \quad \Rightarrow \quad Y \qquad \text{where } P(Y \leq y) = \exp\{-\theta y^{d-\alpha}\}.$$

PROOF. For $y > 0$, compute

$$\nu(\mathbb{R}^d \times [y, \infty)) = \int_{\mathbb{R}^d} \int_y^\infty \nu(dx, dz) = \int_{\mathbb{R}^d} \frac{dx}{(y + q|x|)^\alpha}.$$

Using the substitution $u_1 = x_1 + \cdots + x_d$ and $u_i = x_i$ for $i \geq 2$, and then $y + qu_1 = y/v$, we get

$$\int_{\mathbb{R}^d} \frac{dx}{(y+q|x|)^\alpha} = 2^d \int_0^\infty \frac{u_1^{d-1}}{(y+qu_1)^\alpha} \int_{\substack{u_2+\cdots+u_d \leq u_1 \\ u_i \geq 0}} du_2 \cdots du_n \, du_1$$

$$= \frac{2^d}{(d-1)!} \int_0^\infty \frac{u_1^{d-1} \, du_1}{(y+qu_1)^\alpha}$$

$$= \frac{2^d y^{d-\alpha}}{q^d (d-1)!} \int_0^1 v^{\alpha-d-1}(1-v)^{d-1} \, dv$$

$$= \frac{2^d B(\alpha - d, d)}{y^{\alpha - d} q^d (d-1)!} = \theta y^{d-\alpha}.$$

Since $\nu$ is the intensity measure of the Poisson process $\Pi$, we obtain

$$P(\Pi(\mathbb{R}^d \times [y, \infty)) = 0) = \exp\{-\theta y^{d-\alpha}\}.$$

We eventually have $r(t) > 2dt$ by Lemma 3.1, and hence, by Lemma 2.1(b) that $N(t) = \max_z \Psi_t(z)$ eventually. Using relative compactness of $\mathbb{R}^d \times [y, \infty)$ in $\hat{H}$, we obtain, by Lemma 3.7, that

$$P\left(\frac{N(t)}{a_t} \leq y\right) = P(\Pi_t(\mathbb{R}^d \times [y, \infty)) = 0) \to P(\Pi(\mathbb{R}^d \times [y, \infty)) = 0)$$

$$= \exp\{-\theta y^{d-\alpha}\},$$

which completes the proof. □

3.4. *Almost sure behavior of $N(t)$.* In this section we complete the proof of Theorem 1.1. Taking into account Proposition 3.2, it suffices to prove the following proposition.

PROPOSITION 3.9 [Almost sure bounds on $\log N(t)$].

$$\limsup_{t \to \infty} \frac{\log N(t) - d/(\alpha - d) \log t}{\log \log t} = -\frac{d-1}{\alpha - d} \quad \text{for } d > 1,$$

$$\limsup_{t \to \infty} \frac{\log N(t) - d/(\alpha - d) \log t}{\log \log \log t} = \frac{1}{\alpha - d} \quad \text{for } d = 1$$

and

$$\liminf_{t \to \infty} \frac{\log N(t) - d/(\alpha - d) \log t}{\log \log t} = -\frac{d}{\alpha - d} \quad \text{for } d \geq 1.$$

PROOF. To study the lim sup behavior, we apply, for $t$ large enough, the first estimate of Lemma 3.5 to get, for every $\delta > 0$,

(3.4) $$N(t) \leq \max_{r > 0}\left[r^{d/\alpha}(\log r)^{1/\alpha}(\log \log r)^{1/\alpha}(\log \log \log r)^{1/\alpha + \delta} - \frac{r}{t}\log\frac{r}{2det}\right].$$

Further, we obtain, for every $\delta > 0$,

(3.5) $$N(t) \geq \max_{r > 0}\left[r^{d/\alpha}(\log r)^{1/\alpha}(\log \log r)^{1/\alpha}(\log \log \log r)^{1/\alpha - \delta} - \frac{r}{t}\log\frac{r}{2det}\right],$$

*if* the maximum on the right-hand side is attained at a point $\hat{r}_t$ satisfying the second inequality from Lemma 3.5. This holds for infinitely many $t$ as $\hat{r}_t \to \infty$ continuously in $t$. It therefore remains to prove that the logarithms of the right-hand sides of (3.4) and (3.5) have the required asymptotics.

We deal with both cases simultaneously. For a fixed $\eta$, denote

$$(3.6) \quad f_t(r) = r^{d/\alpha}(\log r)^{1/\alpha}(\log \log r)^{1/\alpha}(\log \log r)^{1/\alpha+\eta} - \frac{r}{t}\log\frac{r}{2det}$$

and denote by $\hat{r}_t$ a maximizer of $f_t$. The condition $d < \alpha$ implies that $\hat{r}_t \to \infty$ so that we have

$$0 = f_t'(\hat{r}_t) = (d/\alpha)\hat{r}_t^{d/\alpha-1}(\log \hat{r}_t)^{1/\alpha}(\log\log \hat{r}_t)^{1/\alpha}(\log\log\log \hat{r}_t)^{1/\alpha+\eta}(1+o(1))$$
$$- \frac{1}{t}\log\frac{\hat{r}_t}{2dt},$$

which, in turn, implies that

$$(3.7) \quad \begin{aligned}(d/\alpha)\hat{r}_t^{d/\alpha}(\log \hat{r}_t)^{1/\alpha}(\log\log \hat{r}_t)^{1/\alpha}(\log\log\log \hat{r}_t)^{1/\alpha+\eta}(1+o(1)) \\ = \frac{\hat{r}_t}{t}\log\frac{\hat{r}_t}{2dt}.\end{aligned}$$

Hence, $\hat{r}_t/t \to \infty$ and taking the logarithm, we obtain $\log(\hat{r}_t/t) = (d/\alpha)[\log \hat{r}_t] \times (1+o(1))$. Substituting this into (3.7) and taking the logarithm, we obtain

$$\log t + \alpha^{-1}\log\log \hat{r}_t + \alpha^{-1}[\log\log\log \hat{r}_t](1+o(1))$$
$$= (1 - d/\alpha)\log \hat{r}_t + \log\log \hat{r}_t.$$

Simplifying, we get

$$(3.8) \quad \begin{aligned}\log t = (1 - d/\alpha)\log \hat{r}_t + (1 - 1/\alpha)\log\log \hat{r}_t \\ - \alpha^{-1}[\log\log\log \hat{r}_t](1+o(1)).\end{aligned}$$

Obviously, we can look for an asymptotic for $\log \hat{r}_t$ in the form

$$\log \hat{r}_t = a_1\log t + a_2\log\log t + a_3[\log\log\log t](1+o(1)).$$

Substituting this into (3.8), we get

$$a_1(1 - d/\alpha) = 1, \qquad a_2(1 - d/\alpha) + 1 - 1/\alpha = 0, \qquad a_3(1 - d/\alpha) - 1/\alpha = 0,$$

which implies

$$(3.9) \quad \log \hat{r}_t = \frac{\alpha}{\alpha - d}\log t - \frac{\alpha - 1}{\alpha - d}\log\log t + \frac{1}{\alpha - d}[\log\log\log t](1+o(1)).$$

Substituting (3.7) into (3.6) and using (3.9), we obtain

$$\log f_t(\hat{r}_t) = \log((1 - d/\alpha)\hat{r}_t^{d/\alpha}(\log \hat{r}_t)^{1/\alpha}(\log\log \hat{r}_t)^{1/\alpha}$$
$$\times (\log\log\log \hat{r}_t)^{1/\alpha+\eta}(1+o(1)))$$

$$= \frac{d}{\alpha-d}\log t + \left[\frac{1}{\alpha} - \frac{d(\alpha-1)}{\alpha(\alpha-d)}\right]\log\log t$$
$$+ \left[\frac{1}{\alpha} + \frac{d}{\alpha(\alpha-d)}\right][\log\log\log t](1+o(1))$$
$$= \frac{d}{\alpha-d}\log t - \frac{d-1}{\alpha-d}\log\log t + \frac{1}{\alpha-d}[\log\log\log t](1+o(1)).$$

Note that the right-hand side does not depend on $\eta$. Therefore, the almost sure upper and lower bounds in (3.4) and (3.5) have identical asymptotics. The second term vanishes for $d = 1$, which explains the difference between the asymptotics for the one-dimensional and multi-dimensional cases.

To study the liminf behavior, we apply the third estimate of Lemma 3.5 and obtain, for $t$ large enough,

$$N(t) \geq \max_{r>0}\left[c_1 r^{d/\alpha}(\log\log r)^{-1/\alpha} - \frac{r}{t}\log\frac{r}{2det}\right].$$

Denote

(3.10) $$f_t(r) = c_1 r^{d/\alpha}(\log\log r)^{-1/\alpha} - \frac{r}{t}\log\frac{r}{2det}$$

and denote by $\hat{r}_t$ a point where the maximum of $r \mapsto f_t(r)$ is achieved. The condition that $d < \alpha$ implies that $\hat{r}_t/t \to \infty$ and we have

$$0 = f_t'(\hat{r}_t) = c_1(d/\alpha)\hat{r}_t^{d/\alpha-1}(\log\log\hat{r}_t)^{-1/\alpha}(1+o(1)) - \frac{1}{t}\log\frac{\hat{r}_t}{2dt},$$

which implies that

(3.11) $$c_1(d/\alpha)\hat{r}_t^{d/\alpha}(\log\log\hat{r}_t)^{-1/\alpha}(1+o(1)) = \frac{\hat{r}_t}{t}\log\frac{\hat{r}_t}{2dt}.$$

Hence, $\hat{r}_t/t \to \infty$, and taking the logarithm, we obtain $\log(\hat{r}_t/t) = (d/\alpha) \times [\log\hat{r}_t](1+o(1))$. Substituting this into (3.11) and taking the logarithm, we obtain

(3.12) $$\log t = (1 - d/\alpha)\log\hat{r}_t + [\log\log\hat{r}_t](1+o(1)).$$

Obviously, we can look for an asymptotic for $\log\hat{r}_t$ in the form

$$\log\hat{r}_t = a_1\log t + a_2[\log\log t](1+o(1)).$$

Substituting this into (3.12), we get $a_1(1 - d/\alpha) = 1$, $a_2(1 - d/\alpha) + 1 = 0$, which implies

(3.13) $$\log\hat{r}_t = \frac{\alpha}{\alpha-d}\log t - \frac{\alpha}{\alpha-d}[\log\log t](1+o(1)).$$

Substituting (3.11) into (3.10) and using (3.13), we obtain

$$\log f_t(\hat{r}_t) = \log\left(c_1\left(1 - \frac{d}{\alpha}\right)\hat{r}_t^{d/\alpha}(\log\log \hat{r}_t)^{-1/\alpha}(1 + o(1))\right)$$
$$= \frac{d}{\alpha - d}\log t - \frac{d}{\alpha - d}[\log\log t](1 + o(1)),$$

which proves that $-d/(\alpha - d)$ is a lower bound for the lim inf. Note that we cannot prove the equality using just estimates for $M_r$, but it follows from the weak convergence proved in Proposition 3.8 as it implies

$$\frac{\log N(t) - d/(\alpha - d)\log t}{\log\log t} \Rightarrow -\frac{d}{\alpha - d}$$

and so there is a sequence $t_n \to \infty$ along which the convergence holds in the almost sure sense. □

**4. Potentials with stretched exponential tails.** We now focus on potentials with distribution function $F(x) = 1 - e^{-x^\gamma}$, $x \geq 0$, for some positive $\gamma < 1$, and include the case $\gamma = 1$ of exponential potentials if this is possible at no additional cost. In Section 4.1 we describe the behavior of the maximum $M_r$ and show that $\underline{N}(t)$ and $N(t)$ are bounds for $L_t$ up to the order $O(1)$. In Section 4.2 we obtain the weak limits theorems for $\underline{N}(t)$ and $N(t)$. In Section 4.3 we show that the approximation of $L_t$ by $\underline{N}(t)$ holds up to order $o((\log t)^{1/\gamma - 1})$. Therefore, the weak limit theorem for $\underline{N}(t)$, which we obtained in Section 4.2, directly implies the weak limit theorem for $L_t$, which completes the proof of Theorem 1.4. Finally, in Section 4.4 we find almost sure estimates for $N(t)$ and $\underline{N}(t)$ using the bounds for $M_r$ and then find the almost sure asymptotics for $L_t$ stated in Theorem 1.3.

4.1. *Approximation of $L_t$ up to constant order.* In this section we first derive bounds for $M_r$ (Lemma 4.1) using the analogous results obtained in the case of Pareto potentials. Then, in Lemma 4.2, we improve the upper and lower bound for $L_t$ obtained in Section 2.1 and show that $L_t$ is squeezed between $\underline{N}(t) - 2d$ and $N(t) - 2d$ up to an additive error which converges to zero.

LEMMA 4.1 (Bounds on $M_r$). *Assume $0 < \gamma \leq 1$ and let $\delta \in (0,1)$ and $c > 0$. Then almost surely, as $r \to \infty$, the following estimates hold:*

$M_r \leq (d\log r)^{1/\gamma} + \gamma^{-1}(d\log r)^{1/\gamma - 1}\log\log r$
$\quad + (\log r)^{1/\gamma - 1}(\log\log r)^\delta \quad$ *eventually for all $r$,*

$M_r \geq (d\log r)^{1/\gamma} + \gamma^{-1}(d\log r)^{1/\gamma - 1}\log\log r \quad$ *for infinitely many $r$,*

$$M_r \geq (d\log r)^{1/\gamma} - (\gamma^{-1} + c)(d\log r)^{1/\gamma-1} \log\log\log r \quad \text{eventually for all } r,$$

$$M_r \leq (d\log r)^{1/\gamma}$$
$$- (\gamma^{-1} - c)(d\log r)^{1/\gamma-1} \log\log\log r \quad \text{for infinitely many } r.$$

*In particular, for $\hat{c} > \gamma^{-1} \log d$, we have $\log M_r \leq \gamma^{-1} \log\log r + \hat{c}$ eventually for all $r$.*

PROOF. This lemma is analogous to Lemma 3.5 for the polynomial potentials, and the proof can be derived from there. Namely, let us pick some $\alpha > d$ and note that $\hat{\xi}(z) = \exp\{\xi(z)^\gamma/\alpha\}$ is a collection of independent Pareto distributed random variables. Denote $\hat{M}_r = \max_{|z| \leq r} \hat{\xi}(z)$. Then

$$M_r = (\alpha \log \hat{M}_r)^{1/\gamma},$$

where the asymptotics for $\hat{M}_r$ are given by Lemma 3.5. We obtain, with $0 < \delta' < 1/\alpha$,

$$M_r \leq [\alpha \log(r^{d/\alpha}(\log r)^{1/\alpha}(\log\log r)^{1/\alpha}(\log\log\log r)^{1/\alpha+\delta'})]^{1/\gamma}$$
$$\leq (d\log r)^{1/\gamma} + \gamma^{-1}(d\log r)^{1/\gamma-1} \log\log r$$
$$+ \gamma^{-1}(1 + \alpha\delta')(d\log r)^{1/\gamma-1} \log\log\log r$$

eventually, which implies the first stated eventual inequality as $\log\log\log t = o((\log\log t)^\delta)$.

Similarly, using the second inequality from Lemma 3.5, we obtain

$$M_r \geq (d\log r)^{1/\gamma} + \gamma^{-1}(d\log r)^{1/\gamma-1} \log\log r$$
$$+ \gamma^{-1}(1 - \alpha\delta')(d\log r)^{1/\gamma-1} \log\log\log r$$

infinitely often, which implies the second stated inequality as $1 - \alpha\delta' > 0$. The third and the fourth inequalities can be proved in the same way. The final statement is an obvious consequence of the more precise first estimate. □

PROPOSITION 4.2 (Bounds on $L_t$). *Assume $0 < \gamma \leq 1$. Then, almost surely,*

$$\underline{N}(t) - 2d + o(1) \leq L_t \leq N(t) - 2d + o(1) \quad \text{eventually for all } t.$$

PROOF. By Lemma 4.1, we have, for any $\eta > 0$, that $M_r \leq r^\eta$ eventually and, hence, the upper bound follows directly from Lemma 2.2. To obtain the lower bound, recall Lemma 2.3 and note that the function $f(\rho) = (1-\rho)\xi(z) - \frac{|z|}{t}\log\frac{|z|}{e\rho t}$ is maximized at $\rho_0 = \frac{|z|}{t\xi(z)}$ and, hence, $\rho = \rho_0$ gives the most powerful estimate in (2.3), unless it drops out of the interval $(0,1)$, whence it violates the condition $|z| < t\xi(z)$. Computing $f(\rho_0) = \xi(z) - \frac{|z|}{t}\log \xi(z)$ gives

$$L_t \geq \max_{|z| < t\xi(z)} \left[\xi(z) - \frac{|z|}{t}\log \xi(z)\right] - 2d + o(1).$$

Observe that for $|z| \geq t\xi(z)$ one has $\xi(z) - \frac{|z|}{t}\log_+ \xi(z) \leq \xi(z) - \xi(z)\log_+ \xi(z) \leq 1$ by properties of the function $x \mapsto x - x\log_+ x$. On the other hand, the maximum taken over all $z$ converges to infinity, as it is equal to $\underline{N}(t)$ by Lemma 2.1(c). Hence,

$$\max_{|z| < t\xi(z)} \left[\xi(z) - \frac{|z|}{t}\log \xi(z)\right] \geq \max_{z \in \mathbb{Z}^d} \left[\xi(z) - \frac{|z|}{t}\log_+ \xi(z)\right] = \underline{N}(t) + o(1). \quad \square$$

4.2. *Weak convergence of $N(t)$ and $\underline{N}(t)$.* In this section we study weak convergence of $N(t)$ and $\underline{N}(t)$ and show that they have the same asymptotics up to the fourth term, which is the first random term in both cases. However, the limiting distribution of the rescaled $N(t)$ and $\underline{N}(t)$ turn out to be different. Our main tool is again the point processes technique developed in Section 2.2. Denote

$$a_t := (d\log t)^{1/\gamma},$$
$$b_t := d(1/\gamma^2 - 1/\gamma)(d\log t)^{1/\gamma - 1}\log \log t,$$
$$c_t := -(d/\gamma)(d\log t)^{1/\gamma - 1}\log \log \log t,$$
$$d_t := (d\log t)^{1/\gamma - 1},$$

which will turn out to be the first four terms in the asymptotics and

$$r_t := t(\log t)^{1/\gamma - 1}(\log \log t)^{-1},$$

which will be the right scaling for the radii $r(t)$ and $\underline{r}(t)$ which maximize $N(t)$ and $\underline{N}(t)$ (in their original definition). Let $G := (-\infty, \infty]$ and define

$$X_{r,z} := \frac{\xi(z) - a_r}{d_r} \quad \text{and} \quad \mu(dy) := \gamma e^{-\gamma y}\,dy,$$

which is a Radon measure on $(-\infty, \infty]$. For $x \in \mathbb{R}$, we have, as $r \to \infty$,

$$r^d P(X_{r,0} \geq x) = r^d P(\xi(0) \geq a_r + xd_r) = r^d \exp\{-(a_r + xd_r)^\gamma\}$$
$$= \exp\{d\log r - d\log r(1 + xd^{-1}(\log r)^{-1})^\gamma\} \longrightarrow e^{-\gamma x}$$
$$= \mu([x, \infty]),$$

and so condition (2.6) is satisfied. Define $\zeta_r$ by (2.7), and, for each $\tau \in \mathbb{R}$ and $q > 0$, define

$$H^q_\tau := \{(x,y) \in \dot{\mathbb{R}}^d \times G : y \geq q|x|/2 + \tau\},$$

which will be our family of state spaces. Denote by $K^q_{\tau,n} := \{(x,y) \in H^q_\tau : |x| \geq n\}$ a family of compact sets in $H^q_\tau$.

LEMMA 4.3 (Convergence of the point processes). *For each $\tau$ and $q$, we have $\zeta_r|_{H^q_\tau} \Rightarrow \zeta|_{H^q_\tau}$, where $\zeta|_{H^q_\tau}$ denotes a Poisson process on $H^q_\tau$ with intensity measure $\eta := \mathrm{Leb}_d \otimes \mu|_{H^q_\tau}$.*

PROOF. It suffices to show that the conditions (i)–(iv) of Lemma 2.4 are satisfied:

(i) As
$$\eta(H^q_\tau) = \int_{\mathbb{R}^d} dx \int_{q|x|/2+\tau}^\infty \gamma e^{-\gamma y}\, dy = e^{-\gamma\tau} \int_{\mathbb{R}^d} e^{-\gamma q|x|/2}\, dx < \infty,$$
the measure $\eta$ is finite on $H^q_\tau$ and hence a Radon measure.

(ii) To prove that $\zeta_r|_{H^q_\tau}$ is a point process in $H^q_\tau$, it suffices to show that, with probability one $\zeta_r(H^q_\tau) < \infty$. It follows from Lemma 4.1 that, for $\delta < 1$, $\xi(z) \leq M_{|z|} < |z|^\delta$ eventually for all $z$, which implies
$$X_{r,z} = \frac{\xi(z) - a_r}{d_r} \leq \frac{|z|^\delta - a_r}{d_r} < \frac{q|z|}{2r} + \tau$$
eventually for all $z$. Hence, $(|z|/r, X_{r,z}) \in H^q_\tau$ for just finitely many $z$.

(iii) Obviously the projection of each compact set in $H^q_\tau$ to the second component is compact in $G$.

(iv) The convergence $\lambda_r \xrightarrow{v} \eta$ with respect to continuous test functions $f : \mathbb{R}^d \times G \to [0, \infty)$ with compact support in $\mathbb{R}^d \times G$ follows from (2.6). To prove the vague convergence on $H^q_\tau$, we need to show that, additionally, for each $\varepsilon > 0$, there is $n$ such that $\lambda_r(K^q_{\tau,n}) < \varepsilon$ eventually for all $r$.

Let $\delta > 0$. Denote $f_r = \varepsilon(\log r)$ for some $\varepsilon > 0$. We have
$$\lambda_r(K^q_{\tau,n}) = \sum_{|z| \geq nr} P\left(X_{r,0} \geq q\frac{|z|}{2r} + \tau\right) = \sum_{|z| \geq nr} P\left(\xi(0) \geq a_r + d_r q\frac{|z|}{2r} + d_r \tau\right)$$
$$= \sum_{|z| \geq nr} \exp\left\{-d\log r\left(1 + \frac{q|z|/(2r) + \tau}{d\log r}\right)^\gamma\right\} =: J_1(n,r) + J_2(r),$$

where $J_1(n,r)$ and $J_2(r)$ correspond to the summations over $nr \leq |z| \leq rf_r$, and $|z| \geq rf_r$, respectively. Using a Taylor expansion, which is valid for $\varepsilon > 0$ sufficiently small, we get

(4.1)
$$J_1(n,r) = e^{-\gamma\tau+o(\varepsilon)} \sum_{nr \leq |z| \leq rf(r)} r^{-d} \exp\left\{-\gamma q\frac{|z|}{2r}\right\}$$
$$\xrightarrow{r \to \infty} e^{-\gamma\tau+o(\varepsilon)} \int_{|x| \geq n} e^{-\gamma q|x|/2}\, dx.$$

Further, for $r$ large enough, we have $|\tau| \leq qf_r/4$ and so $q|z|/(2r) + \tau \geq q|z|/(4r)$ for all $z$ such that $|z| \geq rf_r$. This, together with the monotonicity of $x \mapsto \exp\{-d\log r(1 + \frac{qx}{4rd\log r})^\gamma\}$, implies

$$J_2(r) \leq \sum_{|z| \geq rf_r} \exp\left\{-d\log r\left(1 + \frac{q|z|}{4rd\log r}\right)^\gamma\right\}$$

$$\leq \int_{|x| \geq rf_r - d} \exp\left\{-d\log r\left(1 + \frac{q|x|}{4rd\log r}\right)^\gamma\right\} dx$$

$$= 2\int_{y \geq rf_r - d} y^{d-1} \exp\left\{-d\log r\left(1 + \frac{qy}{4rd\log r}\right)^\gamma\right\} dy$$

$$\leq r^{d+o(1)} \int_{z \geq 1 + q\varepsilon/(8d))(d\log r)^{1/\gamma}} z^{d-1} \exp\{-z^\gamma\} dz,$$

where $z := (d\log r)^{1/\gamma}(1 + qy/(4rd\log r))$, and we have used that

$$y^{d-1} \leq z^{d-1}\{(d\log r)^{1-1/\gamma} 4r/q\}^{d-1} = z^{d-1} r^{d-1+o(1)}.$$

It is not hard to see that the final integral is of order $r^{d(1+q\varepsilon/(8d))+o(1)}$, which proves that $J_2(r) = o(1)$ and, hence, that $J_2(r) < \varepsilon/3$ for $r$ large enough. Since the integral in (4.1) is finite, $n$ can be chosen large enough so that also $J_1(n,r) < \varepsilon/3$ eventually for all $r$. □

Denote

$$(4.2) \quad \Psi_t(z) := \xi(z) - \frac{|z|}{t}\log\frac{|z|}{2det} \quad \text{and} \quad \underline{\Psi}_t(z) := \xi(z) - \frac{|z|}{t}\log_+\xi(z),$$

recalling that $\log_+(x) = \log(x \vee 1)$. Define random variables

$$Y_{t,z} := \frac{\Psi_t(z) - a_{r_t}}{d_{r_t}} \quad \text{and} \quad \underline{Y}_{t,z} := \frac{\underline{\Psi}_t(z) - a_{r_t}}{d_{r_t}}$$

and point processes

$$\Pi_t^\tau := \sum_{\{z \in \mathbb{Z}^d, Y_{t,z} \geq \tau\}} \varepsilon_{(zr_t^{-1}, Y_{t,z})} \quad \text{and} \quad \underline{\Pi}_t^\tau := \sum_{\{z \in \mathbb{Z}^d, \underline{Y}_{t,z} \geq \tau\}} \varepsilon_{(zr_t^{-1}, \underline{Y}_{t,z})}.$$

Finally, denote $q := d^{1-1/\gamma}(\gamma^{-1} - 1)$ and $\underline{q} := d^{1-1/\gamma}\gamma^{-1}$.

LEMMA 4.4 (Convergence of the transformed point processes). *For each $t$, $\Pi_t^\tau$ and $\underline{\Pi}_t^\tau$ are point processes on $\hat{H}_\tau := \dot{\mathbb{R}}^{d+1} \setminus (\mathbb{R}^d \times (-\infty, \tau))$. As $t \to \infty$, we have the following:*

(a) $\Pi_t^\tau \Rightarrow \Pi^\tau$, *where $\Pi^\tau$ is a Poisson process on $\hat{H}_\tau$ with intensity measure $\nu|_{\hat{H}_\tau}$ defined by*

$$\nu(dx, dy) = dx \otimes \gamma e^{-\gamma(y+q|x|)} dy.$$

(b) $\underline{\Pi}_t^\tau \Rightarrow \underline{\Pi}^\tau$, where $\underline{\Pi}^\tau$ is a Poisson process on $\hat{H}_\tau$ with intensity measure $\underline{\nu}|_{\hat{H}_\tau}$ defined by

$$\underline{\nu}(dx, dy) = dx \otimes \gamma e^{-\gamma(y+\underline{q}|x|)} \, dy.$$

PROOF. (a) Observe that

$$Y_{t,z} = \frac{\Psi_t(z) - a_{r_t}}{d_{r_t}} = \frac{\xi(z) - a_{r_t}}{d_{r_t}} - \frac{|zr_t^{-1}|}{td_{r_t}r_t^{-1}} \log \frac{|zr_t^{-1}|}{2detr_t^{-1}}.$$

Since $tr_t^{-1} \to \infty$ and $r_t d_{r_t}^{-1} t^{-1} \log(r_t t^{-1}) \to q$, we can write

$$\Pi_t^\tau = (\zeta_{r_t}|_{H_\tau^q} \circ T_t^{-1})|_{\hat{H}_\tau} \qquad \text{eventually for all } t,$$

for a transformation $T_t : H_\tau^q \to H' := \dot{\mathbb{R}}^{d+1}$ given by

$$T_t : (x, y) \mapsto \begin{cases} (x, y - q|x| - \delta(t, x)), & \text{if } x \neq \infty \text{ and } y \neq \infty, \\ \infty, & \text{otherwise,} \end{cases}$$

where $\delta(t, x) \to 0$ as $t \to \infty$ uniformly on all $(K_{\tau,n}^q)^c$.

We define $T : H_\tau^q \to H'$ by $T(x, y) = (x, y - q|x|)$ if $x \neq \infty$ and $y \neq \infty$ and $T(x, y) = \infty$ otherwise. By Lemmas 4.3 and 2.5, $\zeta_{r_t}|_{H_\tau^q} \circ T_t^{-1} \Rightarrow \zeta|_{H_\tau^q} \circ T^{-1}$. The conditions of the latter are satisfied since (i) is obvious, (ii) is fulfilled as $H_\tau^q$ is itself compact, and (iii) holds as $\delta(x,t) \to 0$ uniformly on $(K_{\tau,n}^q)^c$ and $\eta(K_{\tau,n}^q) \to 0$, recalling that $\eta$ is finite on $H_\tau^q$. Hence, $\Pi_t^\tau$ converges in law to the Poisson process $(\zeta|_{H_\tau^q} \circ T^{-1})|_{\hat{H}_\tau}$ with intensity measure $\eta \circ T^{-1}|_{\hat{H}_\tau} = \nu|_{\hat{H}_\tau}$.

(b) We observe that

$$\underline{Y}_{t,z} = \frac{\Psi_t(z) - a_{r_t}}{d_{r_t}} = \frac{\xi(z) - a_{r_t}}{d_{r_t}} - \frac{|zr_t^{-1}|}{td_{r_t}r_t^{-1}} \log_+ \left\{ a_{r_t} \left[ 1 + \frac{\xi(z) - a_{r_t}}{d_{r_t}} \frac{d_{r_t}}{a_{r_t}} \right] \right\}.$$

Since $d_{r_t}/a_{r_t} \to 0$ and $r_t d_{r_t}^{-1} t^{-1} \log a_{r_t} \to \underline{q}$, we can again represent $\underline{\Pi}_t^\tau$ as a transformation of $\zeta_{r_t}|_{H_\tau^{\underline{q}}}$ via a transformation $T_t$ of the same form with $q$ replaced by $\underline{q}$ and $\delta(t, x)$ by a function $\underline{\delta}(t, x, y)$, which converges uniformly to zero outside the compact sets $\underline{K}_{\tau,n}^q := H_\tau^{\underline{q}} \setminus (\mathbb{R}^d \times (0, n])$. This allows analogous arguments as in (a) to be used in the completion of the proof. □

PROPOSITION 4.5 [Weak asymptotics for $N(t)$ and $\underline{N}(t)$]. *With $\theta$ defined in Theorem 1.4,*

(a) $\frac{N(t) - a_t - b_t - c_t}{d_t} \Rightarrow Y$, *where* $P(Y \leq y) = \exp\{-(1-\gamma)^{-d}\theta e^{-\gamma y}\}$,

(b) $\frac{\underline{N}(t) - a_t - b_t - c_t}{d_t} \Rightarrow \underline{Y}$, *where* $P(\underline{Y} \leq y) = \exp\{-\theta e^{-\gamma y}\}$.

PROOF. We start with the proof of (a). Given any $y \in \mathbb{R}$, we compute

$$\nu(\mathbb{R}^d \times [y, \infty)) = \gamma \int_{\mathbb{R}^d} \int_y^\infty e^{-\gamma(\hat{y}+q|x|)} \, dx \, d\hat{y} = e^{-\gamma y} \int_{\mathbb{R}^d} e^{-\gamma q|x|} \, dx$$

$$= e^{-\gamma y}(2/(\gamma q))^d = (1-\gamma)^{-d}\theta e^{-\gamma y}.$$

Pick $\tau < y$. Since $\nu|_{\hat{H}_\tau}$ is the intensity measure of the Poisson process $\Pi^\tau$, we have

$$P(\Pi^\tau(\mathbb{R}^d \times [0, \infty)) = 0) = \exp\{-\theta e^{-\gamma y}\}.$$

By Lemma 2.1(b), we have $N(t) = \max_z \Psi_t(z)$, provided $r(t) > 2dt$. Let us show that this is fulfilled with probability converging to one. If $r(t) \le 2dt$, then by the monotonicity properties of $M_r$ and $\frac{r}{t} \log \frac{r}{2det}$ one has $r(t) = 2dt$ and $N(t) = M_{2dt} + 2d$. Now we use an eventual almost sure lower bound for $\underline{N}(t)$, and hence for $N(t)$, which will be proved in Lemma 4.10 using only Lemma 2.1 and the bounds provided in Lemma 4.1. By Lemma 4.10, for $0 < c < (1/\gamma^2 - 1/\gamma)d^{1/\gamma}$ and sufficiently large $t$,

$$\begin{aligned}
P(r(t) \le 2dt) &= P(M_{2dt} + 2d \ge \underline{N}(t)) \\
&= P(M_{2dt} \ge (d\log t)^{1/\gamma} + c(\log t)^{1/\gamma - 1}\log\log t) \\
&= 1 - [1 - P(\xi(0) \ge (d\log t)^{1/\gamma} \\
&\qquad + c(\log t)^{1/\gamma - 1}\log\log t)]^{(2dt)^d(1+o(1))\kappa_d}.
\end{aligned}$$
(4.3)

Hence,

$$\begin{aligned}
\log P(r(t) > 2dt) &= -(2dt)^d(1+o(1))\kappa_d \\
&\quad \times \exp\{-[(d\log t)^{1/\gamma} + c(\log t)^{1/\gamma - 1}\log\log t]^\gamma\} \\
&= -\kappa_d \exp\{-d^{1-1/\gamma}c\gamma(1+o(1))\log\log t\} = o(1).
\end{aligned}$$
(4.4)

Now, using relative compactness of $\mathbb{R}^d \times [y, \infty)$ in $\hat{H}_\tau$, we obtain from Lemma 4.4 that

$$\begin{aligned}
P\left(\frac{N(t) - a_{r_t}}{d_{r_t}} \le y, r(t) > 2dt\right) &= P(\Pi^\tau_t(\mathbb{R}^d \times [y, \infty)) = 0, r(t) > 2dt) \\
&\to P(\Pi^\tau(\mathbb{R}^d \times [y, \infty)) = 0) \\
&= \exp\{-(1-\gamma)^{-d}\theta e^{-\gamma y}\}.
\end{aligned}$$

This proves (a) as

$$\begin{aligned}
a_{r_t} &= (d\log t + d(1/\gamma - 1)\log\log t - d\log\log\log t)^{1/\gamma} \\
&= a_t + b_t + c_t + o(d_t),
\end{aligned}$$
(4.5)

$$d_{r_t} = d_t + o(d_t).$$
(4.6)

To prove (b), recall from Lemma 2.1(c) that $\underline{N}(t) = \underline{\Psi}_t(z)$ and argue analogously as in (a). $\square$

4.3. *Sharp approximation of $L_t$ by $\underline{N}(t)$.* In this part we show that up to fourth order, that is, up to $o(d_t)$, we can approximate $L_t$ by $\underline{N}(t)$.

PROPOSITION 4.6 [Sharp approximation of $L_t$ by $\underline{N}(t)$]. *As $t \to \infty$,*
$$P(\underline{N}(t) + o(d_t) \leq L_t \leq \underline{N}(t) + o(d_t)) \to 1.$$

Note that Theorem 1.4 follows directly by combining the weak limit theorem for $\underline{N}(t)$ (Proposition 4.5) and Proposition 4.6. The remainder of this section will be devoted to the proof of Proposition 4.6; we start by introducing some notation.

Denote $M_n^{(0)} := M_n$ and, for $0 < i \leq \kappa_d(n) n^d - 1$, we set
$$M_n^{(i)} := \max\{\xi(z) : |z| \leq n, \xi(z) \neq M_n^{(j)} \ \forall j < \lfloor i \rfloor\}.$$

Let $0 < \rho_1 < \rho_2 < 1$ and
$$k_n := \lfloor n^{\rho_1} \rfloor \quad \text{and} \quad m_n := \lfloor n^{\rho_2} \rfloor.$$

Further, let $f, g : (0, \infty) \to (0, \infty)$ be two monotonically continuous functions going to zero and infinity, respectively. Assume that $\log f_t = o(\log t)$ and $t^{-1} \log g_t = o(d_t)$, that is, that the convergence of $f_t$ to zero and of $g_t$ to infinity is not too fast. Recall that $J_t$ denotes the number of jumps of the random walk $X$ before time $t$. We split $U(t)$ in the form $U(t) = U_1(t) + U_2(t) + U_3(t)$, where

$$U_1(t) := \mathbb{E}_0\bigg[\exp\bigg\{\int_0^t \xi(X_s)\,ds\bigg\} \\ \times \mathbf{1}\bigg\{r_t f_t \leq J_t \leq r_t g_t, \exists i < k_{J_t} \max_{0 \leq s \leq t} \xi(X_s) = M_{J_t}^{(i)}\bigg\}\bigg],$$

$$U_2(t) := \mathbb{E}_0\bigg[\exp\bigg\{\int_0^t \xi(X_s)\,ds\bigg\}\mathbf{1}\bigg\{r_t f_t \leq J_t \leq r_t g_t, \max_{0 \leq s \leq t} \xi(X_s) \leq M_{J_t}^{(k_{J_t})}\bigg\}\bigg],$$

$$U_3(t) := \mathbb{E}_0\bigg[\exp\bigg\{\int_0^t \xi(X_s)\,ds\bigg\}(\mathbf{1}\{J_t < r_t f_t\} + \mathbf{1}\{J_t > r_t g_t\})\bigg].$$

The strategy of the proof of Proposition 4.6 is as follows: We first show that the contributions of $U_2$ and $U_3$ are negligible by showing that, with probability close to one, we have that $U_2 = o(U)$ and $U_3 = o(U)$ (see Lemma 4.8 below). Then we study $U_1$ very carefully and show that, with high probability, $\frac{1}{t} \log U_1(t) \leq \underline{N}(t) + o(d_t)$ (see Lemma 4.9 below). These steps (except studying $U_3$) require detailed information about the asymptotic behavior of the upper order statistics, which is collected in Lemma 4.7.

LEMMA 4.7 (Upper order statistics). *There is a constant $c > 0$ such that, almost surely:*

(i) $M_n^{(0)} - M_n^{(k_n)} \geq c(\log n)^{1/\gamma}$ *eventually for all $n$,*
(ii) $M_n^{(k_n)} - M_n^{(m_n)} \geq c(\log n)^{1/\gamma}$ *eventually for all $n$.*

PROOF. Recall that, by Lemma 4.1, we have
$$\lim_{n \to \infty} M_n^{(0)} (\log n)^{-1/\gamma} = d^{1/\gamma}$$
almost surely. Hence, to prove (i) and (ii), it suffices to show that, for each $\rho \in (0, 1)$
$$\lim_{n \to \infty} M_n^{(n^\rho)} (\log n)^{-1/\gamma} = (d - \rho)^{1/\gamma}.$$

To simplify further computations, note that $\hat\xi(z) = \xi(z)^\gamma$ defines a field of independent exponentially distributed random variables. Then we have $\hat M_n = (M_n)^\gamma$, where $\hat M_n = \max_{|z| \leq n} \hat\xi(z)$ and so it suffices to show that

(4.7) $$\lim_{n \to \infty} \hat M_n^{(n^\rho)} (\log n)^{-1} = d - \rho.$$

Denote by $\ell_n := \kappa_d(n) n^d$ the number of points in the ball of radius $n$. By [5], Proposition 4.1.2, the distribution of $\hat M_n^{(n^\rho)}$ is given by

(4.8) $$P(\hat M_n^{(n^\rho)} \leq x) = \sum_{i=0}^{\lfloor n^\rho \rfloor} \binom{\ell_n}{i} e^{-xi} (1 - e^{-x})^{\ell_n - i}.$$

Then, for each $0 < \varepsilon < d - \rho$, we obtain using $\binom{\ell_n}{i} \leq \ell_n^i$ that

$$P(\hat M_n^{(n^\rho)} \leq (d - \rho - \varepsilon) \log n)$$
$$= \sum_{i=0}^{\lfloor n^\rho \rfloor} \binom{\ell_n}{i} n^{-i(d-\rho-\varepsilon)} (1 - n^{-d+\rho+\varepsilon})^{\ell_n - i}$$
$$\leq (1 - n^{-d+\rho+\varepsilon})^{\ell_n - n^\rho} \sum_{i=0}^{\lfloor n^\rho \rfloor} (\ell_n n^{-d+\rho+\varepsilon})^i$$
$$\leq \exp\{-(1 + o(1)) \ell_n n^{-d+\rho+\varepsilon}\} (n^\rho + 1)(n^{\rho+\varepsilon}(1 + o(1))\kappa_d)^{n^\rho}$$
$$= \exp\{-n^{\rho+\varepsilon}(1 + o(1))\kappa_d\}.$$

Since this sequence is summable, the Borel–Cantelli lemma implies that

(4.9) $$\hat M_n^{(n^\rho)} > (d - \rho - \varepsilon) \log n \quad \text{eventually for all } n.$$

To prove the eventual upper bound, we again use (4.8) and obtain

$$P(\hat{M}_n^{(n^\rho)} \geq (d-\rho+\varepsilon)\log n) = \sum_{\lfloor n^\rho \rfloor + 1}^{\ell_n} \binom{\ell_n}{i} n^{-i(d-\rho+\varepsilon)}(1-n^{-d+\rho+\varepsilon})^{\ell_n - i}$$

$$\leq \sum_{\lfloor n^\rho \rfloor + 1}^{\ell_n} \binom{\ell_n}{i} n^{-i(d-\rho+\varepsilon)}.$$

We can approximate the binomial coefficient using Stirling's formula, which gives

$$\binom{\ell_n}{i} \leq \frac{\ell_n^i}{i!} \leq \left(\frac{e\ell_n}{i}\right)^i$$

for $i$ large enough. This implies, eventually for all $n$,

$$P(\hat{M}_n^{(n^\rho)} \geq (d-\rho+\varepsilon)\log n) \leq \sum_{\lfloor n^\rho \rfloor + 1}^{\ell_n} \left[\frac{e\ell_n}{i n^{d-\rho+\varepsilon}}\right]^i \leq \ell_n \left[\frac{e\ell_n}{n^{d+\varepsilon}}\right]^{n^\rho}$$

$$= e^{-\varepsilon n^\rho (1+o(1))\log n} \leq e^{-n^\rho},$$

where we estimated the sum by its largest term. Again, the sequence is summable and so the Borel–Cantelli lemma implies $\hat{M}_n^{(n^\rho)} < (d-\rho+\varepsilon)\log n$ eventually for all $n$, which, together with the eventual lower bound (4.9), completes the proof of (4.7). □

LEMMA 4.8 [$U_2(t)$ and $U_3(t)$ are negligible]. *As $t \to \infty$,*
$$P(U_2(t) \leq e^{-t}U(t)) \to 1 \quad \text{and} \quad P(U_3(t) \leq e^{-t}U(t)) \to 1.$$

PROOF. First, let us find upper bounds for $U_2$ and $U_3$ similar to the upper bound for $U$ found in Lemma 2.2. More precisely, let us show that almost surely

$$(4.10) \quad \frac{1}{t}\log U_2(t) \leq \max_{\{r_t f_t \leq n \leq r_t g_t\}} \left[M_n^{(k_n)} - \frac{n}{t}\log\frac{n}{2det}\right] - 2d + o(1),$$

$$(4.11) \quad \frac{1}{t}\log U_3(t) \leq \max_{\{n < r_t f_t\} \cup \{n > r_t g_t\}} \left[M_n - \frac{n}{t}\log\frac{n}{2det}\right] - 2d + o(1)$$

eventually for all $t$. The proofs of the two inequalities follow the same line of argument as the proof for the upper bound of $U(t)$. Similarly to (2.1), we have that

$$U_2(t) = \sum_{r_t f_t \leq n \leq r_t g_t} \mathbb{E}_0\left[\exp\left\{\int_0^t \xi(X_s)\,ds\right\}\mathbf{1}\{J_t = n, \max_{0 \leq s \leq t}\xi(X_s) \leq M_n^{(k_n)}\}\right]$$

$$\leq \sum_{r_t f_t \leq n \leq r_t g_t} e^{tM_n^{(k_n)}}\mathbb{P}_0(J_t = n) = \sum_{r_t f_t \leq n \leq r_t g_t} e^{tM_n^{(k_n)} - 2dt}\frac{(2dt)^n}{n!}.$$

Since $M_n^{(k_n)} \leq M_n$, the rest of the computations in Lemma 2.2 is true with $M_n$ replaced by $M_n^{(k_n)}$, which implies the first inequality. In a similar way, we obtain that

$$U_3(t) = \sum_{\{n<r_t f_t\} \cup \{n>r_t g_t\}} \mathbb{E}_0\left[\exp\left\{\int_0^t \xi(X_s)\,ds\right\}\mathbf{1}\{J_t = n\}\right]$$

$$\leq \sum_{\{n<r_t f_t\} \cup \{n>r_t g_t\}} e^{tM_n - 2dt}\frac{(2dt)^n}{n!}.$$

We can again use the computations in Lemma 2.2 with $\mathbb{N}$ replaced by $\{n < r_t f_t\} \cup \{n > r_t g_t\}$ to complete the proof of the second inequality.

Let us now prove the statements of the lemma. Using the lower bound for $U(t)$ from Proposition 4.2, the bound (4.10) and the first estimate from Lemma 4.7, we have eventually

$$\frac{1}{t}\log\frac{U_2(t)}{U(t)} \leq \max_{r_t f_t \leq n \leq r_t g_t}\left[M_n^{(k_n)} - \frac{n}{t}\log\frac{n}{2det}\right] - \underline{N}(t) + o(1)$$

$$\leq \max_{r_t f_t \leq n \leq r_t g_t}\left[M_n - \frac{n}{t}\log\frac{n}{2det}\right] - c(\log(r_t f_t))^{1/\gamma} - \underline{N}(t) + o(1)$$

$$\leq N(t) - \underline{N}(t) - d_t(1 + o(1))cd^{1-1/\gamma}\log t$$

$$= d_t\left(\frac{N(t) - a_t - b_t - c_t}{d_t} - \frac{\underline{N}(t) - a_t - b_t - c_t}{d_t}\right.$$

$$\left. - (1 + o(1))cd^{1-1/\gamma}\log t\right),$$

where the last line converges in law to $-\infty$ by Proposition 4.5. This proves $P(U_2(t) \leq e^{-t}U(t)) \to 1$.

Further, using (4.11) and again the lower bound for $U(t)$ from Proposition 4.2, we have eventually

(4.12) $\quad \frac{1}{t}\log\frac{U_3(t)}{U(t)} \leq \max_{\{n<r_t f_t\} \cup \{n>r_t g_t\}}\left[M_n - \frac{n}{t}\log\frac{n}{2det}\right] - \underline{N}(t) + o(1).$

Denote the maximum in the previous line by $S_t$. Analogously to Lemma 2.1(b), one can prove that

$$S_t \leq \max\left\{M_{2dt} + 2d, \max_{\{|z|<r_t f_t\} \cup \{|z|>r_t g_t\}} \Psi_t(z)\right\},$$

where $\Psi_t(z)$ is defined by (4.2). The computations (4.3) and (4.4) imply that, for any $y > 0$,

(4.13) $\qquad\qquad P(M_{2dt} + 2d \geq a_t + yb_t) = o(1).$

Further, for each $y \in \mathbb{R}$ we pick $\tau < y$ and apply Lemma 4.4 to show that, for each $\varepsilon > 0$, eventually

$$P\left(\max_{\{|z|<r_t f_t\} \cup \{|z|>r_t g_t\}} \Psi_t(z) - a_{r_t} \leq d_{r_t} y\right)$$
(4.14)
$$= P(\Pi_t^\tau(\{|x| < f_t \text{ or } |x| > g_t\} \times [y, \infty)) = 0) \longrightarrow 1,$$

since $f_t \to 0$, $g_t \to \infty$ and $\Pi_t^\tau \to \Pi^\tau$ in law on $\hat{H}_\tau$. Using (4.13), (4.14), and the relations (4.5) and (4.6), we obtain that $(S_t - a_t - b_t - c_t)/d_t \Rightarrow -\infty$. This together with (4.12) implies that

$$\frac{1}{t} \log \frac{U_3(t)}{U(t)} \leq d_t \left(\frac{S_t - a_t - b_t - c_t}{d_t} - \frac{\underline{N}(t) - a_t - b_t - c_t}{d_t}\right) \quad \Rightarrow \quad -\infty,$$

as the second term converges to a nontrivial random variable by Proposition 4.5. □

LEMMA 4.9 [The upper bound on $U_1(t)$]. *As* $t \to \infty$, $P(\frac{1}{t} \log U_1(t) \leq \underline{N}(t) + o(d_t)) \to 1$.

PROOF. In order to obtain such a precise bound, we need to carefully estimate the contribution of different types of trajectories to the Feynman–Kac formula. Denote by

$$\mathcal{P}_n = \{y = (y_0, y_1, \ldots, y_n) : y_0 = 0, |y_{i-1} - y_i| = 1,$$
$$\exists 0 \leq l < k_n \text{ s.t. } \max_{0 \leq i \leq n} \xi(y_i) = M_n^{(l)}\}$$

the set of all discrete time paths in $\mathbb{Z}^d$ with $n$ steps, hitting a point where one of the $k_n$ maximal values of $\xi$ over the ball of radius $n$ is achieved. Let $(\tau_i)_{i \in \mathbb{N}_0}$ be a sequence of independent exponentially distributed random variables with parameter $2d$. Denote by $\mathsf{E}$ the expectation with respect to $(\tau_i)$. Averaging over all random paths following the same geometric path $y$ (with different timings), we obtain

$$U_1(t) = \sum_{r_t f_t \leq n \leq r_t g_t} \sum_{y \in \mathcal{P}_n} (2d)^{-n}$$
(4.15)
$$\times \mathsf{E}\left[\exp\left\{\sum_{i=0}^{n-1} \tau_i \xi(y_i) + \left(t - \sum_{i=0}^{n-1} \tau_i\right) \xi(y_n)\right\}\right.$$
$$\left. \times \mathbf{1}\left\{\sum_{i=0}^{n-1} \tau_i < t, \sum_{i=0}^{n} \tau_i > t\right\}\right].$$

Note that, as $y$ can have self-intersections, some of the values of $\xi$ over $y$ can be the same. We would like to avoid the situation when the maximum of $\xi$ over $y$ is taken at more than one point. Therefore, for each path $y$, we slightly change the potential over $y$. Namely, we denote by $i(y) := \min\{i : \xi(y_i) = \max_{0 \leq j \leq n} \xi(y_j)\}$ the index of the first point where the maximum of the potential over the path is attained. Then we define the modified version of the potential $\xi^y : \{0, \ldots, n\} \to \mathbb{R}$ by

$$\xi_i^y = \begin{cases} \xi(y_i), & \text{if } i \neq i(y), \\ \xi(y_i) + 1, & \text{if } i = i(y). \end{cases}$$

Using $\xi(y_i) \leq \xi_i^y$, we obtain that

$$\mathsf{E}\left[\exp\left\{\sum_{i=0}^{n-1} \tau_i \xi(y_i) + \left(t - \sum_{i=0}^{n-1} \tau_i\right)\xi(y_n)\right\} \mathbf{1}\left\{\sum_{i=0}^{n-1} \tau_i < t, \sum_{i=0}^{n} \tau_i > t\right\}\right]$$

$$\leq \mathsf{E}\left[\exp\left\{\sum_{i=0}^{n-1} \tau_i \xi_i^y + \left(t - \sum_{i=0}^{n-1} \tau_i\right)\xi_n^y\right\} \mathbf{1}\left\{\sum_{i=0}^{n-1} \tau_i < t, \sum_{i=0}^{n} \tau_i > t\right\}\right]$$

$$= (2d)^{n+1} e^{t\xi_n^y} \int_{\mathbb{R}_+^{n+1}} \exp\left\{\sum_{i=0}^{n-1} x_i[\xi_i^y - \xi_n^y - 2d] - 2\, dx_n\right\}$$

(4.16)
$$\times \mathbf{1}\left\{\sum_{i=0}^{n-1} x_i < t, \sum_{i=0}^{n} x_i > t\right\} dx_0 \cdots dx_n$$

$$= (2d)^n e^{-2dt} e^{t\xi_n^y} \int_{\mathbb{R}_+^n} \exp\left\{\sum_{i=0}^{n-1} x_i[\xi_i^y - \xi_n^y]\right\}$$

$$\times \mathbf{1}\left\{\sum_{i=0}^{n-1} x_i < t\right\} dx_0 \cdots dx_{n-1}.$$

The substitution $\hat{x}_i = x_i$ for $i \neq i(y)$ and $\hat{x}_n = t - \sum_{i=0}^{n-1} x_i$ shows that one can interchange the rôle played by $\xi_n^y$ and $\xi_{i(y)}^y$ in the last expression (more generally, one can show that it is symmetric with respect to all indices). We obtain

$$e^{t\xi_n^y} \int_{\mathbb{R}_+^n} \exp\left\{\sum_{i=0}^{n-1} x_i[\xi_i^y - \xi_n^y]\right\} \mathbf{1}\left\{\sum_{i=0}^{n-1} x_i < t\right\} dx_0 \cdots dx_{n-1}$$

$$= e^{t\xi_{i(y)}^y} \int_{\mathbb{R}_+^n} \exp\left\{\sum_{i \neq i(y)} \hat{x}_i[\xi_i^y - \xi_{i(y)}^y]\right\}$$

(4.17)
$$\times \mathbf{1}\left\{\sum_{i \neq i(y)} \hat{x}_i < t\right\} d\hat{x}_0 \cdots d\hat{x}_{i(y)-1}\, d\hat{x}_{i(y)+1}\, d\hat{x}_n$$

$$\leq e^{t\xi_{i(y)}^y} \int_{\mathbb{R}_+^n} \exp\left\{\sum_{i \neq i(y)} \hat{x}_i [\xi_i^y - \xi_{i(y)}^y]\right\} d\hat{x}_0 \cdots d\hat{x}_{i(y)-1}\, d\hat{x}_{i(y)+1}\, d\hat{x}_n$$

$$= e^{t\xi_{i(y)}^y} \prod_{i \neq i(y)} \frac{1}{\xi_{i(y)}^y - \xi_i^y}.$$

Denote by $n(y) = |\{y_0, \ldots, y_n\}|$ the number of different points in the path $y$. Since the path visits $n(y)$ different points, it cannot leave the ball of radius $n(y)$. Therefore, we have that

(4.18) $$\xi_{i(y)}^y \leq M_{n(y)}^{(0)} + 1.$$

By construction, $\xi_{i(y)}^y - \xi_i^y \geq 1$ for all $i \neq i(y)$. Hence, we can drop some terms in the product in order to obtain a further upper bound. Further, since the path visits $n(y)$ different points, there are indices $j_1 > \cdots > j_{n(y)-1}$, which are all different from $i(y)$, such that

$$\xi_{j_i}^y \leq M_{n(y)}^{(i)} \leq M_n^{(i)}, \qquad 1 \leq i < n(y).$$

On the other hand, since $y \in \mathcal{P}_n$, there exists $0 \leq l < k_n$ such that $\xi_{i(y)}^y = M_n^{(l)} + 1$ and so

(4.19) $$\xi_{i(y)}^y - \xi_{j_i}^y \geq M_n^{(l)} + 1 - M_n^{(i)} \geq M_n^{(l)} - M_n^{(i)}, \qquad 1 \leq i < n(y).$$

Note that this estimate becomes trivial for $i \leq l$ and therefore will use it only for $m_n \leq i < n(y)$.

Combining (4.16), (4.17) and using the estimates (4.18) and (4.19), we get

$$\mathsf{E}\left[\exp\left\{\sum_{i=0}^{n-1} \tau_i \xi(y_i) + \left(t - \sum_{i=0}^{n-1} \tau_i\right)\xi(y_n)\right\} \mathbf{1}\left\{\sum_{i=0}^{n-1} \tau_i < t, \sum_{i=0}^{n} \tau_i > t\right\}\right]$$

$$\leq (2d)^n e^{-2dt} e^{t\xi_{i(y)}^y} \prod_{i \neq i(y)} \frac{1}{\xi_{i(y)}^y - \xi_i^y}$$

$$\leq \max_{0 \leq l < k_n} \left[(2d)^n e^{-2dt} e^{t(M_{n(y)}^{(0)}+1)} \prod_{i=m_n}^{n(y)-1} \frac{1}{M_n^{(l)} - M_n^{(i)}}\right]$$

$$\leq (2d)^n e^{-2dt+t} e^{tM_{n(y)}^{(0)}} (M_n^{(k_n)} - M_n^{(m_n)})^{-n(y)+m_n}.$$

Combining this with (4.15), we obtain

$$U_1(t) \leq \sum_{r_t f_t \leq n \leq r_t g_t} \sum_{y \in \mathcal{P}_n} e^{-2dt+t} e^{tM_{n(y)}^{(0)}} (M_n^{(k_n)} - M_n^{(m_n)})^{-n(y)+m_n}$$

$$\leq r_t g_t \max_{r_t f_t \leq n \leq r_t g_t} \max_{1 \leq p \leq n} \exp\{n \log(2d) - 2td + t + tM_p^{(0)}$$
$$- (p - m_n) \log(M_n^{(k_n)} - M_n^{(m_n)})\}.$$

Note that by the first estimate from Lemma 4.1, we have, eventually for all $n$,

$$\max_{1 \leq p \leq n} \log M_p = \log M_n \leq \frac{1}{\gamma} \log \log n + O(1).$$

Using the second statement of Lemma 4.7, we obtain

$$\frac{1}{t} \log U_1(t) \leq \max_{r_t f_t \leq n \leq r_t g_t} \max_{1 \leq p \leq n} \left\{ M_p^{(0)} - \frac{p - m_n}{t} \log(M_n^{(k_n)} - M_n^{(m_n)}) \right\} + o(d_t)$$

$$\leq \max_{r_t f_t \leq n \leq r_t g_t} \max_{1 \leq p \leq n} \left\{ M_p^{(0)} - \frac{p - m_n}{t} \log(c(\log n)^{1/\gamma}) \right\} + o(d_t)$$

$$= \max_{r_t f_t \leq n \leq r_t g_t} \max_{1 \leq p \leq n} \left\{ M_p^{(0)} - \frac{p}{t} \log(\log n)^{1/\gamma} \right\} + o(d_t)$$

$$\leq \max_{r_t f_t \leq n \leq r_t g_t} \max_{1 \leq p \leq n} \left\{ M_p - \frac{p}{t} \log M_p \right\} + o(d_t)$$

$$\leq \max_{p > 0} \left\{ M_p - \frac{p}{t} \log M_p \right\} + o(d_t) = \underline{N}(t) + o(d_t),$$

which completes the proof. □

PROOF OF PROPOSITION 4.6. The lower bound has been proved in Proposition 4.2. To prove the upper bound, recall that $U(t) = U_1(t) + U_2(t) + U_3(t)$ and so

$$L_t = \frac{1}{t} \log U_1(t) + \frac{1}{t} \log\left(1 + \frac{U_2(t)}{U_1(t)} + \frac{U_3(t)}{U_1(t)}\right).$$

By Lemma 4.8, the last term converges to zero in law, and, by Lemma 4.9, the probability that the first term is bounded by $\underline{N}(t) + o(d_t)$ converges to 1. □

4.4. *Almost sure asymptotics.* In this section we prove Theorem 1.3. In order to do so, we find bounds for $\underline{N}(t)$ and $N(t)$. First, we find a lower bound for $\underline{N}(t)$ (and so on $L_t$) which holds eventually and observe that it coincides with the first two terms of the weak asymptotic of $L_t$ proved in Theorem 1.4. This gives us the liminf asymptotic for $L_t$. Further, we find a lower bound for $\underline{N}(t)$ which holds infinitely often and an eventual upper bound for $N(t)$, and see that the bounds coincide. This gives the limsup asymptotic for $L_t$ as it is squeezed between $\underline{N}(t)$ and $N(t)$.

The next three lemmas provide estimates for the bounds of $N(t)$ and $\underline{N}(t)$. The proofs all follow the same pattern: one replaces $M_r$ in the definition of $N(t)$ or $\underline{N}(t)$ by its upper or lower bound computed in Lemma 4.1, and then finds the maximum of the new deterministic function. Note that this approach cannot be used to find an upper bound for $N(t)$, which would hold infinitely often.

LEMMA 4.10 [Eventual lower bound on $\underline{N}(t)$]. *Let $\gamma \leq 1$. Then*

$$\underline{N}(t) \geq (d \log t)^{1/\gamma} + (1/\gamma^2 - 1/\gamma)d^{1/\gamma}(\log t)^{1/\gamma - 1} \log \log t$$
$$+ o((\log t)^{1/\gamma - 1} \log \log t),$$

*eventually for all $t$.*

PROOF. By Lemma 2.1(a), the maximum $\underline{N}(t)$ is attained at $\underline{r}(t)$, which goes to infinity. Hence, we can use Lemma 4.1, specifically the third estimate for $M_r$ and the last statement for $\log M_r$, to obtain

$$\underline{N}(t) \geq \max_{r \geq 0} \bigg[ (d \log r)^{1/\gamma} - (\gamma^{-1} + c)(d \log r)^{1/\gamma - 1} \log \log \log r$$
$$- \frac{r}{\gamma t} \log \log r - \frac{r\hat{c}}{t} \bigg].$$

Denoting by $f_t(r)$ the expression in the square brackets, we have

$$(4.20) \quad f'_t(r) = \frac{d^{1/\gamma}(\log r)^{1/\gamma - 1}(1 + o(1))}{\gamma r} - \frac{[\log \log r](1 + o(1))}{\gamma t},$$

where $o(1)$ here is a function of $r$. Denoting by $\hat{r}_t$ a maximizer of $f_t$, we have $\hat{r}_t \to \infty$. $f'_t(\hat{r}_t) = 0$ implies $\hat{r}_t = t\varphi(\hat{r}_t)$, where

$$(4.21) \quad \varphi(r) = d^{1/\gamma}(\log r)^{1/\gamma - 1}(\log \log r)^{-1}(1 + o(1)).$$

This implies that

$$(4.22) \quad \log \varphi(r) = (1/\gamma - 1) \log \log r + o(\log \log r).$$

Using $\hat{r}_t = t\varphi(\hat{r}_t)$, we get $\log \hat{r}_t = \log t + \log \varphi(\hat{r}_t) = \log t + O(\log \log \hat{r}_t)$, which implies

$$(4.23) \quad \log t / \log \hat{r}_t = 1 + o(1),$$

and this yields $\log \varphi(\hat{r}_t) / \log t = o(1)$. Finally, using (4.21), (4.22), (4.23) and considering only the terms up to order $(\log t)^{1/\gamma - 1} \log \log t$, we get

$$f(\hat{r}_t) = (d \log(t\varphi(\hat{r}_t)))^{1/\gamma} - \gamma^{-1}\varphi(\hat{r}_t) \log \log(\hat{r}_t) - \varphi(\hat{r}_t) \log c$$
$$+ o((\log t)^{1/\gamma - 1} \log \log t)$$

$$= (d \log t)^{1/\gamma} + \gamma^{-1} d^{1/\gamma} (\log t)^{1/\gamma - 1} \log \varphi(\hat{r}_t) + o((\log t)^{1/\gamma - 1} \log \log t)$$

$$= (d \log t)^{1/\gamma} + (1/\gamma^2 - 1/\gamma) d^{1/\gamma} (\log t)^{1/\gamma - 1} \log \log t$$

$$\quad + o((\log t)^{1/\gamma - 1} \log \log t),$$

which completes the proof. $\square$

LEMMA 4.11 [i.o. lower bound on $\underline{N}(t)$]. *Let $\gamma \leq 1$. Then*

$$\underline{N}(t) \geq (d \log t)^{1/\gamma} + [(1/\gamma) d^{1/\gamma - 1} + (1/\gamma^2 - 1/\gamma) d^{1/\gamma}](\log t)^{1/\gamma - 1}$$

$$\quad \times [\log \log t](1 + o(1))$$

*infinitely often.*

PROOF. Similarly to the proof of Lemma 4.10, by Lemma 2.1(a), we can use the estimates from Lemma 4.1 for $M_r$ in the definition of $\underline{N}(t)$. Using the second estimate for $M_r$ itself and the last estimate for $\log M_r$, we obtain

$$\underline{N}(t) \geq \max_{r \geq 0} \left[ (d \log r)^{1/\gamma} + \gamma^{-1} (d \log r)^{1/\gamma - 1} \log \log r - \frac{r}{\gamma t} \log \log r - \frac{r \hat{c}}{t} \right],$$

if the maximum of the expression in the square brackets [which we denote by $f_t(r)$] is attained at a point $\hat{r}_t$ such that

(4.24) $\qquad M_{\hat{r}_t} \geq (d \log \hat{r}_t)^{1/\gamma} + \gamma^{-1} (d \log \hat{r}_t)^{1/\gamma - 1} \log \log \hat{r}_t.$

Note that $f'_t(r)$ has the same form (4.20) as in Lemma 4.10, and $\hat{r}_t \to \infty$. Therefore, we can use the same computation as in Lemma 4.10 and show that $\hat{r}_t = t \varphi(\hat{r}_t)$, where $\varphi$ satisfies (4.21), and hence, (4.22) and (4.23) are also satisfied. Using them, we get, considering only the terms up to order $(\log t)^{1/\gamma - 1} \log \log t$,

$$f(\hat{r}_t) = (d \log(t \varphi(\hat{r}_t)))^{1/\gamma} + \gamma^{-1} (d \log(t \varphi(\hat{r}_t)))^{1/\gamma - 1} \log \log(t \varphi(\hat{r}_t))$$

$$\quad + o((\log t)^{1/\gamma - 1} \log \log t)$$

$$= (d \log t)^{1/\gamma} + \gamma^{-1} d^{1/\gamma} (\log t)^{1/\gamma - 1} \log \varphi(\hat{r}_t)$$

$$\quad + \gamma^{-1} d^{1/\gamma - 1} (\log t)^{1/\gamma - 1} [\log \log t](1 + o(1))$$

$$= (d \log t)^{1/\gamma}$$

$$\quad + [(1/\gamma) d^{1/\gamma - 1} + (1/\gamma^2 - 1/\gamma) d^{1/\gamma}](\log t)^{1/\gamma - 1} [\log \log t](1 + o(1)).$$

It remains to check that condition (4.24) holds infinitely often. This is true by Lemma 4.1, using that $t \mapsto \hat{r}_t$ is a continuous function converging to infinity. $\square$

LEMMA 4.12 [Eventual upper bound on $N(t)$]. *Let $\gamma \leq 1$. Then*

$$N(t) \leq (d\log t)^{1/\gamma} + [(1/\gamma)d^{1/\gamma-1} + (1/\gamma^2 - 1/\gamma)d^{1/\gamma}](\log t)^{1/\gamma-1}$$
$$\times [\log\log t](1 + o(1))$$

*eventually for all $t$.*

PROOF. As before, by Lemma 2.1(a), we can use the estimates from Lemma 4.1 for $M_r$ in the definition of $N(t)$. Using the first estimate, we obtain

$$N(t) \leq \max_{r>0}\left[(d\log r)^{1/\gamma} + \gamma^{-1}(d\log r)^{1/\gamma-1}\log\log r\right.$$
$$\left. + (\log r)^{1/\gamma-1}(\log\log r)^{\delta} - \frac{r}{t}\log\frac{r}{2det}\right].$$

Denote by $r \mapsto f_t(r)$ the function in the square brackets and denote by $\hat{r}_t$ a point where $f_t$ attains its maximum. Note that $\hat{r}_t \to \infty$ and compute

$$f'_t(r) = \frac{d^{1/\gamma}(\log r)^{1/\gamma-1}(1+o(1))}{\gamma r} - \frac{1}{t}\log\frac{r}{2dt}.$$

Since $f'_t(\hat{r}_t) = 0$, we obtain $\hat{r}_t = t\varphi(t, \hat{r}_t)$, where

$$\varphi(t, \hat{r}_t) = \gamma^{-1}d^{1/\gamma}(\log \hat{r}_t)^{1/\gamma-1}\left(\log\frac{\hat{r}_t}{2dt}\right)^{-1}(1+o(1))$$

as $\hat{r}_t \to \infty$. Using this, we get

$$\frac{\hat{r}_t}{t}\log\frac{\hat{r}_t}{2dt} = \gamma^{-1}d^{1/\gamma}(\log \hat{r}_t)^{1/\gamma-1}(1+o(1))$$

and, hence, $\varphi(t, \hat{r}_t) = \hat{r}_t/t \to \infty$ for $\gamma < 1$ and is bounded for $\gamma = 1$. Taking the logarithm of the last equality, we obtain $\log\varphi(t, \hat{r}_t) + \log\log\varphi(t, \hat{r}_t) = (1/\gamma - 1)\log\log \hat{r}_t + o(\log\log \hat{r}_t)$. Dividing by $\log \hat{r}_t$, we get $\log\varphi(t, \hat{r}_t) = o(\log \hat{r}_t)$, which implies that $\log\varphi(t, \hat{r}_t) = o(\log t)$ since we have from $\hat{r}_t = t\varphi(t, \hat{r}_t)$ that $\log \hat{r}_t = \log t + \log\varphi(t, \hat{r}_t)$. Finally, we obtain, considering only the terms up to order $(\log t)^{1/\gamma-1}\log\log t$, that

$$f_t(\hat{r}_t) = (d\log(t\varphi(t, \hat{r}_t)))^{1/\gamma} + \gamma^{-1}(d\log(t\varphi(t, \hat{r}_t)))^{1/\gamma-1}\log\log(t\varphi(t, \hat{r}_t))$$
$$+ o((\log t)^{1/\gamma-1}\log\log t)$$
$$= (d\log t)^{1/\gamma} + \gamma^{-1}d^{1/\gamma}(\log t)^{1/\gamma-1}\log\varphi(t, \hat{r}_t)$$
$$+ \gamma^{-1}d^{1/\gamma-1}(\log t)^{1/\gamma-1}[\log\log t](1+o(1))$$
$$= (d\log t)^{1/\gamma} + [(1/\gamma)d^{1/\gamma-1} + (1/\gamma^2 - 1/\gamma)d^{1/\gamma}]$$
$$\times (\log t)^{1/\gamma-1}[\log\log t](1+o(1)),$$

which completes the proof. □

We are finally ready to describe the almost sure behavior of $L_t$ up to second order.

PROOF OF THEOREM 1.3. We know from Proposition 4.2 that $\underline{N}(t) + O(1) \leq L_t \leq N(t) + O(1)$ eventually for all $t$. The first statement of the theorem now follows from Lemmas 4.11 and 4.12. To prove the second one, note that Theorem 1.4 implies

$$\frac{L_t - (d \log t)^{1/\gamma}}{(d \log t)^{1/\gamma - 1} \log \log t} \quad \Rightarrow \quad d(1/\gamma^2 - 1/\gamma),$$

and so $d(1/\gamma^2 - 1/\gamma)$ is an upper bound for the liminf. Since $L_t \geq \underline{N}(t) + O(1)$, the equality follows now from Lemma 4.10. □

**Acknowledgments.** We would like to thank Marek Biskup for his suggestion to use point process techniques, and Hubert Lacoin for providing the proof of Lemma 4.7 and numerous further suggestions on how to shorten our proofs.

R. VAN DER HOFSTAD
EINDHOVEN UNIVERSITY OF TECHNOLOGY
DEPARTMENT OF MATHEMATICS
  AND COMPUTER SCIENCE
P.O. BOX 513
5600 MB EINDHOVEN
THE NETHERLANDS
E-MAIL: rhofstad@win.tue.nl

P. MÖRTERS
UNIVERSITY OF BATH
DEPARTMENT OF MATHEMATICAL SCIENCES
BATH BA2 7AY
UNITED KINGDOM
E-MAIL: maspm@bath.ac.uk

N. SIDOROVA
UNIVERSITY COLLEGE LONDON
DEPARTMENT OF MATHEMATICS
GOWER STREET
LONDON WC1E 6BT
UNITED KINGDOM
E-MAIL: n.sidorova@ucl.ac.uk